\documentclass[journal]{IEEEtran}
\usepackage{amsmath,epsfig,amsfonts,amssymb}
\usepackage{graphics,psfrag,theorem,calc,subfigure,url,bm,cite,color}
\usepackage{amssymb}
\usepackage{amsmath}
\usepackage{latexsym}
\usepackage{mathrsfs}
\usepackage{algorithm}
\usepackage{algorithmic}
\usepackage{bm}
\usepackage{color}
\usepackage{verbatim}
\usepackage{url}
\usepackage{graphicx}  
\usepackage{stfloats}  
\usepackage{tcolorbox}
\tcbuselibrary{skins, breakable, theorems}
\usepackage{pgfplots}
\usepackage{tikz}

\usepackage{hyperref}
\hypersetup{
	colorlinks=true,
	linkcolor=blue,
	filecolor=magenta,      
	urlcolor=cyan,
}

\newenvironment{proof}{{\bf Proof\,\,}}{\endproof\par}

\newcounter{spb}
\def \openbox{$\sqcup\llap{$\sqcap$}$}
\def \endproof{\enskip \null \nobreak \hfill \openbox \par}

\newtheorem{definition}{Definition}

\newtheorem{defi}{Definition}
\newtheorem{prop}{Proposition}

\newtheorem{fact}{Fact}

\newcommand{\limto}{\rightarrow}
\newcommand{\del}{\partial}

\newcommand{\R}{\mathbb R}

\ifCLASSINFOpdf
\else
\fi

\begin{document}
%
\title{Understanding Notions of Stationarity in Non-Smooth Optimization}
\author{
	Jiajin Li,
	Anthony Man-Cho So, \IEEEmembership{Senior Member, IEEE}, and
	Wing-Kin Ma, \IEEEmembership{Fellow, IEEE}		
	
	\thanks{This work is supported in part by the Hong Kong Research Grants Council (RGC) General Research Fund (GRF) projects CUHK 14208117 and CUHK 14208819, and in part by the CUHK Research Sustainability of Major RGC Funding Schemes project 3133236.}

    \thanks{J. Li and A. M.-C. So is with the Department of Systems Engineering and Engineering Management, The Chinese University of Hong Kong, Shatin, N.T., Hong Kong. E-mail: \{jjli, manchoso\}@se.cuhk.edu.hk.}
	
    \thanks{W.-K. Ma is with the Department of Electronic Engineering, The Chinese University of Hong Kong, Shatin, N.T., Hong Kong. E-mail: wkma@ieee.org.}
}

\markboth{}%
{Shell \MakeLowercase{\textit{et al.}}: Bare Demo of IEEEtran.cls for IEEE Journals}

\maketitle

\begin{abstract}
Many contemporary applications in signal processing and machine learning give rise to structured non-convex non-smooth optimization problems that can often be tackled by simple iterative methods quite effectively. One of the keys to understanding such a phenomenon---and, in fact, one of the very difficult conundrums even for experts---lie in the study of ``stationary points'' of the problem in question. Unlike smooth optimization, for which the definition of a stationary point is rather standard, there is a myriad of definitions of stationarity in non-smooth optimization. In this article, we give an introduction to different stationarity concepts for several important classes of non-convex non-smooth functions and discuss the geometric interpretations and further clarify the relationship among these different concepts. We then demonstrate the relevance of these constructions in some representative applications and how they could affect the performance of iterative methods for tackling these applications.
\end{abstract}

\begin{IEEEkeywords}
non-smooth analysis, subdifferential, stationarity.
\end{IEEEkeywords}

%
\IEEEpeerreviewmaketitle

\section{Introduction}
In recent years, we have witnessed a fast-growing body of literature that utilizes non-convex non-smooth optimization techniques to tackle machine learning and signal processing applications. 
Although such a development seems to run contrary to the long-held belief that non-convex optimization problems pose serious analytic and algorithmic challenges, it is proven to be practically relevant and opens up an exciting avenue for dealing with contemporary applications. For instance, various low-rank matrix recovery problems admit natural non-convex optimization formulations that can be readily tackled by lightweight first-order methods (e.g., (sub)gradient descent or block coordinate descent) and are more scalable than their convex approximations; see, e.g.,~\cite{KBV09,GLM16,TBS+16,LCZL19,LZSV20}. On the other hand, many modern statistical estimation problems involve non-convex loss functions and/or regularizers. While such problems are non-convex, they possess certain convexity properties (which can be made precise) that can be exploited in computation, and there are algorithms that can compute solutions to these problems with good empirical performance; see, e.g.,~\cite{LW12,LW15,APX17,CPS18,NPR19}. Another example that has drawn immense interest is deep neural networks with non-smooth activation functions (e.g., the rectified linear unit (ReLU) $x\mapsto \max\{x,0\}$). To train such networks, one often needs to optimize a loss function that is recursively defined via compositions of linear mappings with nonlinear activation functions. Despite the non-convexity and possible non-smoothness of the loss function, various stochastic algorithms (e.g., stochastic (sub)gradient descent or Adam-type algorithms) for optimizing it can still yield exciting empirical performance on a host of machine learning tasks; see, e.g.,~\cite{KSH12,KB15,CLSH19}. There are many other applications whose natural optimization formulations are non-convex yet highly structured, such as dictionary learning~\cite{SQW17a,SQW17b}, non-negative matrix factorization~\cite{Lin07,LLR16}, and phase retrieval~\cite{BEB18,ELB18}. It is becoming increasingly clear that by carefully exploiting the structure of the non-convex formulation at hand, one can design algorithms that have better empirical performance and runtime than those for solving the corresponding convex approximations.

To better understand such phenomenon, a general approach is to study the ``stationary points'' of the problem in question and investigate how existing iterative methods behave around these stationary points. For smooth optimization, the definition of a stationary point is rather standard. Indeed, consider the unconstrained minimization problem
\begin{equation} \label{eq:opt-prob}
\inf_{\bm{x}\in\R^n} f(\bm{x})
\end{equation}
with $f:\R^n\limto\R$. Suppose that $f$ is smooth and let $\nabla f:\R^n\limto\R^n$ be its gradient. A point $\bm{x}\in\R^n$ is said to be stationary if
\[ \nabla f(\bm{x}) = \bm{0}, \]
which means that $\bm{x}$ is either a local minimum, a local maximum, or a saddle point. However, for non-smooth optimization, one can find a myriad of definitions of a stationary point in the literature; see, e.g.,~\cite{RW04,M13a} and the references therein. It is far from clear how these different definitions of stationarity are related and, more fundamentally, why they need to be introduced. This not only creates potential confusion among readers but also obscures the nature of the solutions that are being computed by different iterative methods.

In this paper, our main objective is to give an introduction to the theory of subdifferentiation for non-convex non-smooth functions, with a focus on motivating the different constructions of the subdifferential and developing the corresponding stationarity concepts for several important function classes, as well as discussing the geometric interpretations and further clarifying the relationship among the different constructions. We will also demonstrate the relevance of these constructions in some representative applications and how they could affect the performance of iterative methods for tackling these applications. Readers may just be as intrigued by what classes of iterative algorithms can lead to efficient computation of a stationary point under the aforementioned concepts. Unfortunately, owing to the need for exposition of more sophisticated concepts and also to the page limitation, we decide not to cover algorithms in this introductory article.

\section{Convex Non-smooth Functions} \label{sec:cvx}
To set the stage for our later developments, let us review the theory of subdifferentiation for convex non-smooth functions.
For simplicity, we restrict our discussion to finite-valued convex functions $f:\R^n\limto\R$. Recall that if $f$ is convex and smooth, then its gradient $\nabla f$ at $\bm{x}\in\R^n$ provides an affine minorant of $f$ at $\bm{x}\in\R^n$; i.e., 
\[ f(\bm{y}) \ge f(\bm{x}) + \nabla f(\bm{x})^T(\bm{y}-\bm{x}) \quad \mbox{for all } \bm{y} \in \R^n. \]
In the non-smooth case, a suitable generalization of gradient is the notion of \emph{subgradient}; i.e., a vector $\bm{s}\in\R^n$ is a subgradient of $f$ at $\bm{x}\in\R^n$ if
\[ f(\bm{y}) \ge f(\bm{x}) + \bm{s}^T(\bm{y}-\bm{x}) \quad \mbox{for all } \bm{y} \in \R^n. \]
Since the subgradient at a point may not be unique, we are led to the notion of \emph{subdifferential}, which is the set
\begin{align}
\del f(\bm{x}) &= \left\{ \bm{s} \in \R^n: f(\bm{y}) \ge f(\bm{x}) + \bm{s}^T(\bm{y}-\bm{x}) \right. \nonumber \\
&\qquad\qquad\qquad \left. \mbox{for all } \bm{y} \in \R^n \right\}.  \label{eq:subdiff-2}
\end{align}

As it turns out, the subdifferential~\eqref{eq:subdiff-2} can be constructed by considering the directional derivative of $f$.
Given a point $\bm{x}\in\R^n$ and a direction $\bm{d}\in\R^n$, the \emph{difference quotient} $q$ of $f$ at $\bm{x}$ is defined by
\begin{equation}\label{eq:diffquo}
t\mapsto q(t) = \frac{f(\bm{x}+t\bm{d})-f(\bm{x})}{t} \quad\mbox{for}~ t>0. 
\end{equation}
Observe that by the convexity of $f$, the function $q$ is increasing in $t$ (see, e.g.,~\cite[Chapter 0, Proposition 6.1]{HL01}) and bounded around 0 (see, e.g.,~\cite[Chapter B, Theorem 3.1.2]{HL01}). Thus, the \emph{directional derivative} of $f$ at $\bm{x}\in\R^n$ in the direction $\bm{d}\in\R^n$, which is defined by
\begin{equation} \label{eq:dir-dev}
f'(\bm{x},\bm{d}) = \lim_{t\searrow0} \frac{f(\bm{x}+t\bm{d})-f(\bm{x})}{t},
\end{equation}
exists and is equal to $f'(\bm{x},\bm{d}) = \inf_{t>0} q(t)$. One of the key properties of $f'$ is the following:
\begin{fact}
	(\cite[Chapter D, Proposition 1.1.2]{HL01}) For any $\bm{x}\in\R^n$, the function $\bm{d}\mapsto f'(\bm{x},\bm{d})$ is finite sublinear (recall that a function $h:\R^n\limto\R\cup\{+\infty\}$ is sublinear if it is convex and satisfies $h(t\bm{x})=t \cdot h(\bm{x})$ for all $\bm{x}\in\R^n$ and $t>0$).
\end{fact}
A fundamental result in convex analysis is that there is a correspondence between closed sublinear functions and closed convex sets; see~\cite[Chapter C]{HL01}. In particular, upon invoking~\cite[Chapter C, Theorem 3.1.1]{HL01}, we know that $f'(\bm{x},\cdot)$ is the \emph{support function} of the non-empty closed convex set
\begin{equation} \label{eq:subdiff-1}
\del f(\bm{x}) = \left\{ \bm{s} \in \R^n : \bm{s}^T\bm{d} \le f'(\bm{x},\bm{d}) \mbox{ for all } \bm{d} \in \R^n \right\};
\end{equation}
i.e., $f'(\bm{x},\bm{d}) = \sup_{\bm{s}\in\del f(\bm{x})} \bm{s}^T\bm{d}$. 
By the finiteness of $f'(\bm{x},\cdot)$, the set $\del f(\bm{x})$ is bounded (\cite[Chapter C, Proposition 2.1.3]{HL01}). Hence, $\del f(\bm{x})$ is in fact compact. It can be shown that~\eqref{eq:subdiff-1} and~\eqref{eq:subdiff-2} describe the same set; see~\cite[Chapter D, Theorem 1.2.2]{HL01}. Interestingly, even though we define the set~\eqref{eq:subdiff-2} without reference to differentiation, its support function turns out to be the directional derivative $f'(\bm{x},\cdot)$.

In applications we often need to compute an element of the subdifferential of a given function. Let us now give the subdifferentials of some concrete convex functions $f$.
\begin{itemize}
\item[--] {\bf (Smooth function).} Suppose that $f$ is differentiable at $\bm{x}$. Then, $\partial f(\bm{x}) = \{\nabla f(\bm{x})\}$; see~\cite[Chapter D, Corollary 2.1.4]{HL01}.

\item[--] {\bf (Norm).} Let $f$ be a norm on $\R^n$. Then, 
\[ \del f(\bm{x}) = \{ \bm{s} \in\R^n: \bm{s}^T\bm{x} = f(\bm{x}), \, f_*(\bm{s}) \le 1\}, \]
where $f_*$ is the dual norm of $f$ defined by $f_*(\bm{s}) = \sup_{\bm{d}\in\R^n:f(\bm{d})\le1} \bm{d}^T\bm{s}$; see~\cite[Chapter D, Example 3.1]{HL01}.

\smallskip
In particular, for the $\ell_1$-norm $f(\cdot)=\|\cdot\|_1$, we have $\del(\|\bm{x}\|_1) = {\rm Sign}(\bm{x})$, where ${\rm Sign}$ is the element-wise sign function given by
\[ [{\rm Sign}(\bm{x})]_i = \left\{
\begin{array}{c@{\quad}l}
\{x_i/|x_i|\} & \mbox{if } x_i \not= 0, \\
\noalign{\smallskip}
[-1,1] & \mbox{otherwise};
\end{array}
\right.
\]
for the $\ell_2$-norm $f(\cdot)=\|\cdot\|_2$, we have 
\[ \del(\|\bm{x}\|_2) = \left\{
\begin{array}{c@{\quad}l}
\{ \bm{x}/\|\bm{x}\|_2 \} & \mbox{if } \bm{x}\not=\bm{0}, \\
B(\bm{0},1) & \mbox{otherwise},
\end{array}
\right.
\]
where $B(\bm{0},1)$ is the unit ball centered at the origin.

\item[--] {\bf (Max function).} Suppose that $f$ takes the form $f(\cdot) = \max_{\bm{y}\in Y} g(\cdot,\bm{y})$, where $Y\subseteq\R^\ell$ is compact and $g:\R^n\times Y \limto\R$ is such that $\R^n\ni\bm{x}\mapsto g(\bm{x},\bm{y})$ is convex for each $\bm{y}\in Y$ and $Y \ni \bm{y}\mapsto g(\bm{x},\bm{y})$ is continuous for each $\bm{x}\in\R^n$. Let $Y(\bm{x}) = \{ \bm{y} \in Y : f(\bm{x}) = g(\bm{x},\bm{y}) \}$ be the set of optimal solutions to $\max_{\bm{y} \in Y} g(\bm{x},\bm{y})$. Then,
\begin{equation} \label{eq:max-sd}
\del f(\bm{x}) = {\rm conv}\left\{ \bigcup_{\bm{y} \in Y(\bm{x})} \del g(\bm{x},\bm{y}) \right\};
\end{equation}
cf.~\cite[Chapter D, Theorem 4.4.2]{HL01}.

\smallskip
The above result is extremely useful, as many convex functions can be represented as the maximum of a collection of convex functions. For instance, let $\mathcal{S}^n$ denote the set of $n\times n$ real symmetric matrices and consider the largest eigenvalue function $\mathcal{S}^n \ni \bm{M} \mapsto \lambda(\bm{M})$. By the Courant-Fischer theorem, we have the characterization
\[ \lambda(\bm{M}) = \max_{\bm{u}\in\R^n:\|\bm{u}\|_2=1} \bm{u}^T\bm{M}\bm{u}. \]
Since the function $\bm{M}\mapsto \bm{u}^T\bm{M}\bm{u}$ is linear with gradient $\bm{u}\bm{u}^T$, it follows from~\eqref{eq:max-sd} that
\[ \del\lambda(\bm{M}) = {\rm conv}\left\{ \bm{u}\bm{u}^T: \|\bm{u}\|_2=1, \, \bm{M}\bm{u} = \lambda(\bm{M})\bm{u} \right\}. \]

\item[--] {\bf (Sum rule).} Suppose that $f$ takes the form $f=\alpha_1f_1+\alpha_2f_2$, where $f_1,f_2:\R^n\limto\R$ are convex functions and $\alpha_1,\alpha_2>0$ are positive scalars. Then, $\del f=\alpha_1\del f_1 + \alpha_2\del f_2$; see~\cite[Chapter D, Theorem 4.1.1]{HL01}.

\item[--] {\bf (Composition with affine mapping).} Suppose that $f$ takes the form $f = g\circ A$, where $g:\R^m\limto\R$ is a convex function and $A:\R^n\limto\R^m$ is an affine mapping given by $A(\bm{x}) = \bm{A}_0\bm{x}+\bm{b}$ with $\bm{A}_0\in\R^{m\times n}$ and $\bm{b}\in\R^m$. Then,
\[ \del f(\bm{x}) = \bm{A}_0^T\del g(A(\bm{x})) = \left\{ \bm{A}_0^T\bm{s} : \bm{s} \in \del g(A(\bm{x})) \right\}; \]
see~\cite[Chapter D, Theorem 4.2.1]{HL01}. The above result can be viewed as a \emph{chain rule} for subdifferentials. Note that we restrict ourselves to the composition of a convex function with an affine mapping here, as the resulting function is guaranteed to be convex and hence its subdifferential~\eqref{eq:subdiff-2} is well defined. To obtain more general chain rules, we need to define a notion of subdifferential for non-convex functions. This will be our objective in subsequent sections.

\item[--] {\bf (Indicator).} Although our development so far focuses on finite-valued convex functions, it can be extended to convex functions that take values in $\R\cup\{+\infty\}$. One important example of such functions is the indicator of a closed convex set. Specifically, let $C\subseteq\R^n$ be a closed convex set and define the indicator of $C$ by
\begin{equation} \label{eq:ind}
\mathbb{I}_C(\bm{x}) = \left\{
\begin{array}{c@{\quad}l}
0 & \mbox{if } \bm{x} \in C, \\
+\infty & \mbox{otherwise}.
\end{array}
\right.
\end{equation}
Using the construction~\eqref{eq:subdiff-2} of the subdifferential, it can be verified that 
\begin{equation} \label{eq:ind-sd}
\del\mathbb{I}_C(\bm{x}) = \{ \bm{s} \in \R^n: \bm{s}^T(\bm{y}-\bm{x}) \le 0 \mbox{ for all } \bm{y}\in C\}
\end{equation}
if $\bm{x}\in C$ and $\del\mathbb{I}_C(\bm{x}) = \emptyset$ otherwise. The set on the right-hand side of~\eqref{eq:ind-sd} is known as the \emph{normal cone} to $C$ at $\bm{x}$ and is denoted by $\mathcal{N}_C(\bm{x})$. Each element $\bm{s}\in\mathcal{N}_C(\bm{x})$ is called a \emph{normal direction} to $C$ at $\bm{x}$. The terminology is motivated by the observation that for every $\bm{s}\in\mathcal{N}_C(\bm{x})$, the set $C$ is completely contained in the halfspace $\{ \bm{y}\in\R^n: \bm{s}^T(\bm{y}-\bm{x})\le0 \}$, whose boundary is a hyperplane that passes through $\bm{x}$ and has normal $\bm{s}$; see the figure below.

\vspace{2mm}
\begin{figure}[h]
\centering
\tikzset{every picture/.style={line width=0.75pt}} 

\begin{tikzpicture}[x=0.75pt,y=0.75pt,yscale=-1,xscale=1]

\draw  [fill={rgb, 255:red, 155; green, 155; blue, 155 }  ,fill opacity=0.79 ] (158.9,139.62) .. controls (151.48,116.77) and (154.84,95.2) .. (166.4,91.45) .. controls (177.95,87.7) and (193.33,103.19) .. (200.75,126.04) .. controls (208.17,148.89) and (204.81,170.45) .. (193.25,174.21) .. controls (181.69,177.96) and (166.31,162.47) .. (158.9,139.62) -- cycle ;
\draw    (105.5,160.2) -- (156.91,140.65) ;
\draw  [fill={rgb, 255:red, 0; green, 0; blue, 0 }  ,fill opacity=1 ] (156.9,139.65) .. controls (156.88,138.54) and (157.77,137.64) .. (158.87,137.62) .. controls (159.98,137.61) and (160.88,138.49) .. (160.9,139.59) .. controls (160.91,140.7) and (160.03,141.61) .. (158.93,141.62) .. controls (157.82,141.64) and (156.91,140.75) .. (156.9,139.65) -- cycle ;

\draw  [fill={rgb, 255:red, 155; green, 155; blue, 155 }  ,fill opacity=1 ] (363.76,140.88) -- (297.13,131.13) -- (348,87) ;
\draw  [fill={rgb, 255:red, 0; green, 0; blue, 0 }  ,fill opacity=1 ] (294.9,131.65) .. controls (294.88,130.54) and (295.77,129.64) .. (296.87,129.62) .. controls (297.98,129.61) and (298.88,130.49) .. (298.9,131.59) .. controls (298.91,132.7) and (298.03,133.61) .. (296.93,133.62) .. controls (295.82,133.64) and (294.91,132.75) .. (294.9,131.65) -- cycle ;
\draw  [fill={rgb, 255:red, 155; green, 155; blue, 155 }  ,fill opacity=0.22 ][line width=0.75]  (259,90) -- (295.73,131.05) -- (288.41,185.65) ;

\draw (257,128) node [anchor=north west][inner sep=0.75pt]  [font=\scriptsize]  {$\mathcal{N}_{C}( \bm{x})$};
\draw (304.75,123.59) node [anchor=north west][inner sep=0.75pt]  [font=\footnotesize,rotate=-359.22]  {$\bm{x}$};
\draw (333.46,106.43) node [anchor=north west][inner sep=0.75pt]  [font=\footnotesize,rotate=-359.22]  {$C$};
\draw (99,164.2) node [anchor=north west][inner sep=0.75pt]  [font=\scriptsize]  {$\mathcal{N}_{C}(\bm{x})$};
\draw (162.75,134.00) node [anchor=north west][inner sep=0.75pt]  [font=\footnotesize,rotate=-359.22]  {$\bm{x}$};
\draw (172.46,106.43) node [anchor=north west][inner sep=0.75pt]  [font=\footnotesize,rotate=-359.22]  {$C$};
\end{tikzpicture}
	\caption{Normal cone of a closed convex set.}
\end{figure}
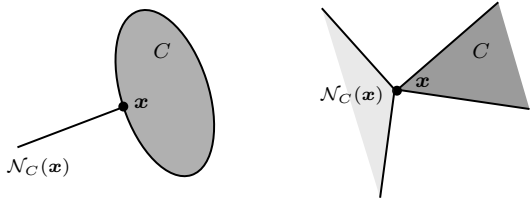
\end{itemize}
\vspace{2mm}


Using the notion of subdifferential, we can formulate the optimality condition of the minimization of a convex non-smooth function. Specifically, let $g:\R^n\limto\R$ be a convex function and $C\subseteq\R^n$ be a closed convex set. Consider the problem
\begin{equation} \label{eq:con-opt}
\inf_{\bm{x}\in C} g(\bm{x}), 
\end{equation}
which can be put into the form~\eqref{eq:opt-prob} by letting $f=g+\mathbb{I}_C$. We then have the following result:
\begin{fact} (cf. {\cite[Theorem 8.15]{RW04}})  \label{fact:cvx-opt}
The following are equivalent:
\begin{enumerate}
\item[(a)] $\bar{\bm x}$ is an optimal solution to~\eqref{eq:con-opt}.

\item[(b)] $\bm{0} \in \del g(\bar{\bm x}) + \mathcal{N}_C(\bar{\bm x})$ (cf.~\eqref{eq:ind-sd}).

\item[(c)] $g'(\bar{\bm x},\bm{y}-\bar{\bm x}) \ge 0$ for all $\bm{y}\in C$.
\end{enumerate}
\end{fact}

Although the main focus of this paper is on notions of stationarity, let us briefly digress and discuss the algorithmic aspects of Problem~\eqref{eq:con-opt}. We say that $\bm{d}\in\R^n$ is a \emph{descent direction} of the convex function $g$ at $\bar{\bm x}\in\R^n$ if there exists a $\bar{t}>0$ satisfying $g(\bar{\bm x} + t \bm{d}) < g(\bar{\bm x})$ for all $t \in (0,\bar{t})$. As can be easily verified, this is equivalent to $g'(\bar{\bm x},\bm{d})<0$. In view of Fact~\ref{fact:cvx-opt}(c), we are thus motivated to use \emph{feasible descent} methods to solve Problem~\eqref{eq:con-opt}. Roughly speaking, at the current iterate $\bm{x}^k\in C$, such methods find a direction $\bm{d}^k$ and step size $\alpha_k>0$ such that the next iterate $\bm{x}^{k+1} = \bm{x}^k + \alpha_k\bm{d}^k $ satisfies $g(\bm{x}^{k+1}) < g(\bm{x}^k)$ and $\bm{x}^{k+1}\in C$. As simple as the above description may seem, there are various subtleties in its implementation. For instance, since $g'({\bm x}^k,\bm{d})<0$ is equivalent to $\max_{\bm{s}\in\del g({\bm x}^k)} \bm{s}^T\bm{d}<0$ (see~\eqref{eq:subdiff-1}), one may be tempted to compute the entire subdifferential $\del g({\bm x}^k)$ in each iteration. However, this could be rather expensive. Moreover, even for the unconstrained minimization of a convex non-smooth function, some natural descent methods (such as a straightforward extension of the steepest descent method for smooth minimization) are not necessarily convergent; see, e.g.,~\cite{Wolfe75}. It turns out that the above difficulties can be overcome. We refer the reader to~\cite{Kiwiel85} for developments in this direction.

\begin{figure}[H]
	\centering
	\includegraphics[width=0.4\textwidth]{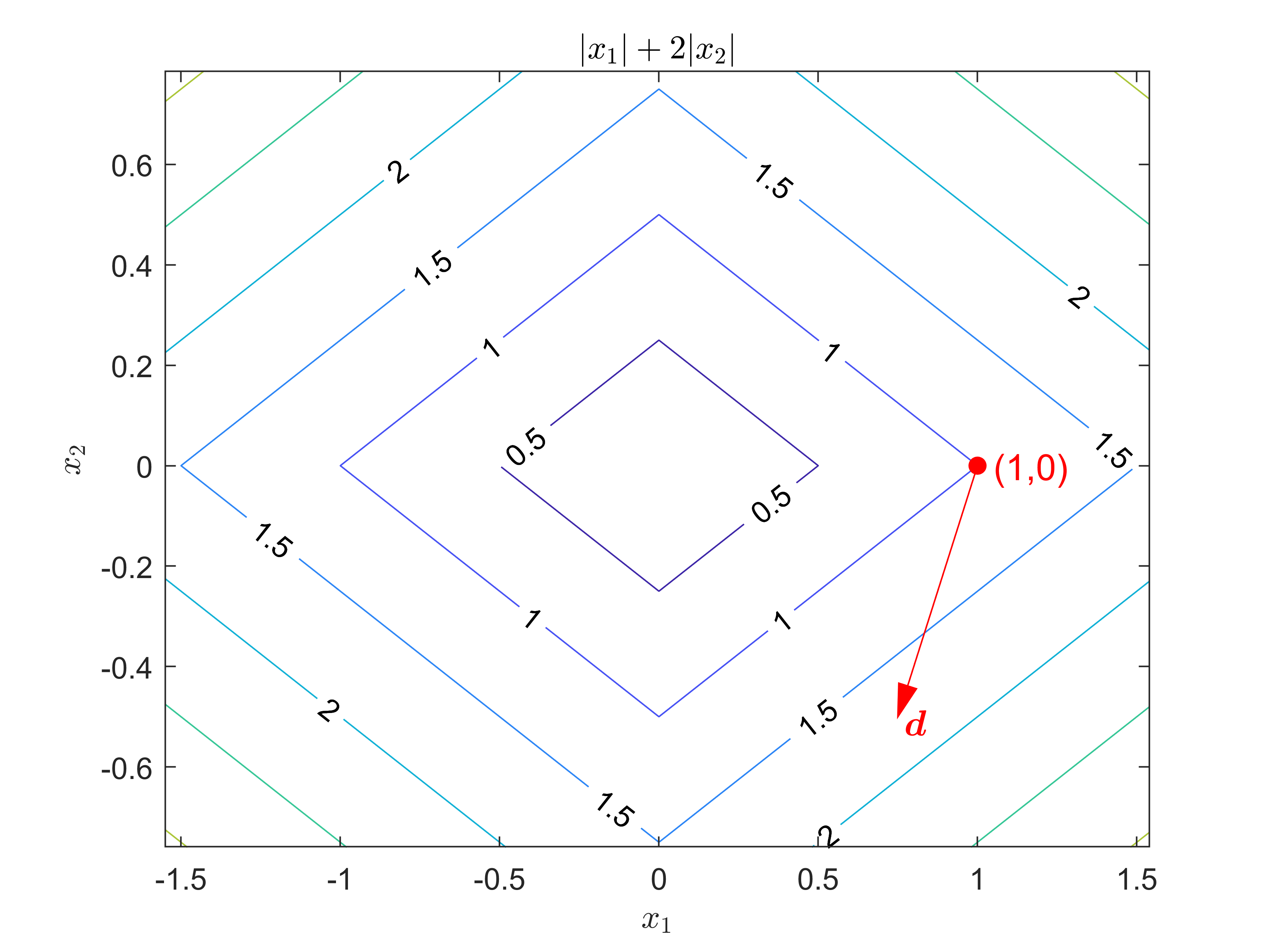}	
	\caption{Contour plot of $f(x_1,x_2) = |x_1| + 2|x_2|$.}
	\label{fig:sub}
\end{figure}

Another idea for solving Problem~\eqref{eq:con-opt} is to use projected subgradient methods. At the current iterate $\bm{x}^k\in C$, such methods proceed by first finding a subgradient $\bm{s}^k \in \del g(\bm{x}^k)$ and choosing a step size $\alpha_k>0$, and then obtaining the next iterate via $\bm{x}^{k+1} = \Pi_C(\bm{x}^k - \alpha_k\bm{s}^k)$, where $\Pi_C$ is the projector onto $C$. It should be noted that subgradient methods are generally \emph{not} descent methods. For instance, consider the function $\R^2 \ni (x_1,x_2) \mapsto f(x_1,x_2) = |x_1| + 2|x_2|$, whose contour plot is given in Figure~\ref{fig:sub}. It can be easily seen that $(1,2) \in \del f(1,0)$, but $\bm{d} = -(1,2)$ is not a descent direction.

\noindent In spite of this, by choosing step sizes that decay at an appropriate rate, it can be shown that subgradient methods will converge to an optimal solution and their convergence rates can be estimated. We refer the reader to~\cite{Goffin77,Shor85} for details.

So far we have only discussed constructions of the subdifferential for convex functions. It should not take long for one to realize that those constructions do not yield much useful information when applied to even some very simple non-convex functions. For instance, if we consider the smooth non-convex function $\R\ni x\mapsto f(x) = -x^2$, then using the definition~\eqref{eq:subdiff-2} we have $\del f(0)=\emptyset$. Another example is the non-smooth non-convex function $\R\ni x\mapsto f(x) = -|x|$, where $\del f(0) = \{s\in\R: s\ge1 ~\mbox{and}~ s\le-1 \} = \emptyset$ according to~\eqref{eq:subdiff-2}. In view of these examples, it is natural to ask whether one can construct a subdifferential that can better capture the geometry of non-smooth non-convex functions. Before we address this question, let us list some desirable properties that we wish such a \emph{generalized subdifferential} to possess. 
First, the subdifferential should be a singleton consisting of the gradient (resp.~coincide with the usual convex subdifferential) when the function in question is smooth (resp.~convex). Second, from a computational point of view, the subdifferential should satisfy some basic calculus rules, particularly the chain rule for composite functions and the sum rule for sum of functions. Without such rules, many concrete non-convex functions that arise in applications cannot be tackled easily. Third, the subdifferential should yield a necessary condition for local optimality; i.e., if $f$ attains a local minimum at $\bar{\bm x}$, then $\bm{0} \in \del f(\bar{\bm x})$. Fourth, the subdifferential should be \emph{tight}, in the sense that the set $\{\bm{x}\in\R^n: \bm{0} \in \del f(\bm{x})\}$ of \emph{stationary points} of $f$ should contain as few non-local minima as possible.
In summary, we have the following desiderata of a generalized subdifferential:
\vspace{1mm}
\begin{tcolorbox}[title = {Desirable Properties of a Generalized Subdifferential}]
	\label{box:1}
\begin{itemize}
	\item[--] for smooth $f$, $\partial f(\bm{x}) = \{ \nabla f(\bm{x})\}$
	\item[--] coincide with the usual convex subdifferential for convex functions
	\item[--] basic calculus rules
	\begin{itemize}
		\item {\it chain rule}: for $f = g \circ F$, $\partial f (\bm{x}) = (JF(\bm{x}))^T \partial g(F(\bm{x}))$, where $JF$ is the Jacobian of $F$ (see~\eqref{eq:jac} for the definition)
		\item {\it sum rule}: for $f = f_1 + f_2$, $\partial f(\bm{x}) = \partial f_1(\bm{x})+\partial f_2(\bm{x})$
	\end{itemize}
	\item[--] necessary condition for local optimality
	\item[--] tight subdifferential
\end{itemize}
\end{tcolorbox}
\vspace{1mm}

The above discussion suggests that one can consider an axiomatic approach to constructing subdifferentials with the desired properties for more general functions. Such an approach has been explored, e.g., in~\cite{Ioffe12}. Another approach, which is more geometric in nature and follows our development for convex functions, is to construct a convex set that serves as the generalized subdifferential and take its support function to be the generalized directional derivative. Alternatively, one can define a sublinear function that serves as the generalized directional derivative and take the convex set it supports as the generalized subdifferential. Let us now take this geometric approach as the starting point of our exposition.


\section{Locally Lipschitz Functions} \label{sec:loc-lip}
As we move beyond convex functions, one direction to explore is the class of locally Lipschitz functions. Such a class captures a wide variety of non-convex functions~\cite{Cla75} and includes the class of convex functions as a special case~\cite[Theorem 10.4]{R97}. Let us recall the definition:
\begin{definition} \label{def:loclip} 
A function $f:\R^n\limto\R$ is locally Lipschitz if for any bounded $S\subseteq\R^n$, there exists a constant $L>0$ such that
	\[ |f(\bm{x})-f(\bm{y})| \le L \|\bm{x}-\bm{y}\|_2 \quad\mbox{for all } \bm{x}, \bm{y} \in S. \]
\end{definition}
By a classic result of Rademacher, a locally Lipschitz function $f$ is differentiable almost everywhere (a.e.)~\cite[Theorem 9.60]{RW04}. In particular, every neighborhood of $\bm{x}$ contains a point $\bm{y}$ for which $\nabla f(\bm{y})$ exists, so that there is at least one cluster point due to the Lipschitzian property. This motivates the following construction, which is known as the \emph{Bouligand subdifferential}:
\[ 
\del_B f({\bm x}) = \left\{ 
\bm{s} \in \R^n:
\begin{aligned}
&\exists \bm{x}^k\limto \bm{x}, \, \nabla f(\bm{x}^k) \mbox{ exists}, \\
&\nabla f(\bm{x}^k) \limto \bm{s} 
\end{aligned}
\right\}.
\]
As a quick illustration, consider the absolute value function $\R\ni x\mapsto f(x)=|x|$. It can be easily verified that $\del_B f(0) = \{-1,1\}$. Such an example is instructive, as it highlights two drawbacks of the Bouligand subdifferential. First, the Bouligand subdifferential does not coincide with the usual convex subdifferential when the function in question is convex. Second, the condition $\bm{0}\in\del_B f(\bm{x})$ is not even necessary for the local optimality of ${\bm x}$. One possible remedy is to \emph{convexify} the Bouligand subdifferential by considering its convex hull; i.e.,
\begin{equation} \label{eq:subdiff-chull}
\del_C f(\bm{x}) =  {\rm conv}(\del_B f({\bm x})).
\end{equation}
It can be shown that $\del_C f(\bm{x})$ so defined is a non-empty compact convex set and is called the \emph{Clarke subdifferential} in the literature; see \cite[Definition 1.1]{Cla75}. From our discussion in Section~\ref{sec:cvx}, we know that $\del_C f(\bm{x})$ can also be described by its support function. This leads to the question: What is the support function of $\del_C f(\bm{x})$ when $f$ is locally Lipschitz? The following remarkable result due to Clarke furnishes the answer.
\begin{fact} \label{fact:clarke}
	(cf.~\cite[Proposition 1.4]{Cla75}) Given a point $\bm{x}\in\R^n$ and a direction $\bm{d}\in\R^n$, the \emph{Clarke directional derivative} of $f$ at $\bm{x}$ in the direction $\bm{d}$ is defined by
	\begin{equation} \label{eq:cla-dir-dev}
	\begin{aligned}
	f^\circ(\bm{x},\bm{d}) &= \limsup_{\bm{x}'\limto \bm{x}, \, t \searrow 0} \frac{f(\bm{x}'+t\bm{d})-f(\bm{x}')}{t} \\
	&= \inf_{\epsilon>0,\atop\lambda>0} \sup_{\bm{x}'\in \bm{x}+\epsilon B(\bm{0},1), \atop t\in(0,\lambda)}  \frac{f(\bm{x}'+t\bm{d})-f(\bm{x}')}{t}.
	\end{aligned}
	\end{equation}
	Then, $f^\circ(\bm{x},\cdot)$ is the support function of the set $\del_C f(\bm{x})$ defined in~\eqref{eq:subdiff-chull}; i.e.,
	\[ f^\circ(\bm{x},\bm{d}) = \max_{\bm{s}\in\del_C f(\bm{x})} \bm{s}^T\bm{d}. \]
	In particular, we have
	\begin{equation} \label{eq:cla-subdiff}
	\del_C f(\bm{x}) = \left\{ \bm{s} \in \R^n : \bm{s}^T\bm{d} \le f^\circ(\bm{x},\bm{d}) \mbox{ for all } \bm{d} \in \R^n \right\}, 
	\end{equation}
	and the function $ \bm{d}\mapsto f^\circ(\bm{x},\bm{d})$ is finite sublinear for all $\bm{x}\in\R^n$. Additionally, we have
	$ f'(\bm{x},\bm{d}) \leq f^\circ(\bm{x},\bm{d})$ if  the directional derivative of $f$ exists.
\end{fact}

We remark that a locally Lipschitz function may not be directionally differentiable. In other words, the difference quotient in \eqref{eq:diffquo} may not have a limit even though it is bounded due to the Lipschitzian property. Here, we give an example to showcase such possibility.
\vspace{2mm}
\begin{tcolorbox}
 Consider the function
\[ \R \ni x \mapsto f(x) =\left\{
\begin{array}{c@{\quad}l}
 x \sin (\log (\frac{1}{x}))& \mbox{if } x >0, \\
0 & {\mbox{otherwise}}.
\end{array} \right.
\]
It is clear that $f$ is smooth on $\R\setminus \{0\}$. Its derivative at any $x>0$ is given by $f'(x) =\sin (\log (\frac{1}{x}))-\cos (\log (\frac{1}{x}))$, which is bounded by $2$. Using this and the structure of $f$, it can be shown that $f$ is locally Lipschitz. However, the directional derivative of $f$ at $\bar{x}=0$ does not exist. In fact, the difference quotient $q(t) = \frac{f(t)}{t} = \sin(\log(\frac{1}{t}))$ does not converge, as can be seen by considering the sequence $t_n = e^{-(n+\frac{1}{2}) \pi}$ and computing
\[q(t_n) = \sin\left( \left( n+\frac{1}{2} \right) \pi \right) =  \left\{
\begin{array}{c@{\quad}l}
1 & \mbox{if $n$ is even}, \\
-1 & {\mbox{otherwise.}}
\end{array} \right.\]
\end{tcolorbox}

It is instructive to compare the two notions of directional derivatives in~\eqref{eq:dir-dev} and~\eqref{eq:cla-dir-dev} from a geometric point of view. The former considers the variation of $f$ along a ray emanating from $\bm{x}$ in the direction $\bm{d}$ (i.e., $f(\bm{x}+t_k\bm{d})$ vs. $f(\bm{x})$ with $t_k\searrow0$), while the latter considers the variation of $f$ in the direction $\bm{d}$ for points in the neighborhood of $\bm{x}$ (i.e., $f(\bm{x}^k+t_k\bm{d})$ vs. $f(\bm{x}^k)$ with $t_k\searrow0$ and $\bm{x}^k\limto \bm{x}$). In particular, the latter is able to explore the behavior of $f$ in a neighborhood of $\bm{x}$ rather than just along a ray emanating from $\bm{x}$. Generally, $ f^\circ(\bm{x},\bm{d})$ is an upper bound on the difference quotient in the neighborhood of $\bm{x}$.   As we shall see, such an idea turns out to be very fruitful when studying the local behavior of non-smooth functions.

Our discussion above reveals a fundamental difference in the theory of subdifferentiation for convex functions and non-convex functions. Specifically, in the convex case, subdifferentiation entails linearization of the function at hand; in the non-convex case, subdifferentiation can be seen as a convexification process. This allows the use of concepts from convex analysis to study the subdifferentials of non-convex functions.

Recall that in Section~\ref{sec:cvx}, we have introudced several properties that the generalized subdifferential should possess. Now, let us check whether the Clarke subdifferential possesses those properties.
\begin{itemize}
	\item[--] {\bf (Smooth function).} If $f$ is \emph{smooth} (i.e., continuously differentiable) at $\bm{x}$, then $f^\circ(\bm{x},\bm{d}) = \nabla f(\bm{x})^T\bm{d}$ for all $\bm{d}\in\R^n$ and $\del_C f(\bm{x}) = \{ \nabla f(\bm{x})\}$; see~\cite[Proposition 1.13]{Cla75}. 

	\item[--] {\bf (Convex function).} As mentioned above, convex functions are locally Lipschitz. In this case, the Clarke subdifferential and Clarke directional derivative take on particularly simple forms. Indeed, the Clarke subdifferential coincides with the usual convex subdifferential \eqref{eq:subdiff-2} due to~\cite[Theorems 17.2 and 25.6]{R97}. In addition, the directional derivative of a convex function, which always exists, is equal to the Clarke directional derivative; i.e., 
	\begin{equation} \label{eq:cvx-cla}
	f^\circ(\bm{x},\bm{d}) = f'(\bm{x},\bm{d}).
	\end{equation}
	
	\item[--] {\bf (Sum rule).} The following example demonstrates that the sum rule $\del_C(f_1+f_2) = \del_Cf_1 + \del_Cf_2$ does not hold in general. Consider the function $f:\R\limto\R$ given by $f(x) = \max\{x,0\} + \min\{0,x\}$.  Let us compute $\del_C f_1(0)$, $\del_Cf_2(0)$, and $\del_Cf(0)$:
	\vspace{3mm}
	\begin{tcolorbox}
	\tikzset{every picture/.style={line width=0.3pt}} 
	\begin{tikzpicture}[x=0.35pt,y=0.35pt,yscale=-0.9,xscale=0.9]
	
	\draw [color={rgb, 255:red, 73; green, 135; blue, 206 }  ,draw opacity=1 ][line width=1.5]  (31,128.58) -- (209.83,128.58)(119.83,48) -- (119.83,203.58) (202.83,123.58) -- (209.83,128.58) -- (202.83,133.58) (114.83,55) -- (119.83,48) -- (124.83,55) (150.83,123.58) -- (150.83,133.58)(181.83,123.58) -- (181.83,133.58)(88.83,123.58) -- (88.83,133.58)(57.83,123.58) -- (57.83,133.58)(114.83,97.58) -- (124.83,97.58)(114.83,66.58) -- (124.83,66.58)(114.83,159.58) -- (124.83,159.58)(114.83,190.58) -- (124.83,190.58) ;
	\draw   ;
	\draw [color={rgb, 255:red, 208; green, 2; blue, 27 }  ,draw opacity=1 ][line width=1.5]    (34.83,128.58) -- (119.83,128.58);
	\draw [color={rgb, 255:red, 208; green, 2; blue, 27 }  ,draw opacity=1 ][line width=1.5]    (119.83,128.58) -- (179.83,64.58);
	\draw [color={rgb, 255:red, 73; green, 135; blue, 206 }  ,draw opacity=1 ][line width=1.5]  (237,128.58) -- (415.83,128.58)(325.83,48) -- (325.83,203.58) (408.83,123.58) -- (415.83,128.58) -- (408.83,133.58) (320.83,55) -- (325.83,48) -- (330.83,55) (356.83,123.58) -- (356.83,133.58)(387.83,123.58) -- (387.83,133.58)(294.83,123.58) -- (294.83,133.58)(263.83,123.58) -- (263.83,133.58)(320.83,97.58) -- (330.83,97.58)(320.83,66.58) -- (330.83,66.58)(320.83,159.58) -- (330.83,159.58)(320.83,190.58) -- (330.83,190.58) ;
	\draw   ;
	\draw [color={rgb, 255:red, 208; green, 2; blue, 27 }  ,draw opacity=1 ][line width=1.5]    (262.83,189.25) -- (325.83,128.58) ;
	\draw [color={rgb, 255:red, 208; green, 2; blue, 27 }  ,draw opacity=1 ][line width=1.5]    (325.83,128.58) -- (389.83,128.25) ;
	\draw [color={rgb, 255:red, 73; green, 135; blue, 206 }  ,draw opacity=1 ][line width=1.5]  (440,130.58) -- (618.83,130.58)(528.83,50) -- (528.83,205.58) (611.83,125.58) -- (618.83,130.58) -- (611.83,135.58) (523.83,57) -- (528.83,50) -- (533.83,57) (559.83,125.58) -- (559.83,135.58)(590.83,125.58) -- (590.83,135.58)(497.83,125.58) -- (497.83,135.58)(466.83,125.58) -- (466.83,135.58)(523.83,99.58) -- (533.83,99.58)(523.83,68.58) -- (533.83,68.58)(523.83,161.58) -- (533.83,161.58)(523.83,192.58) -- (533.83,192.58) ;
	\draw   ;
	\draw [color={rgb, 255:red, 208; green, 2; blue, 27 }  ,draw opacity=1 ][line width=1.5]    (471.83,189.25) -- (528.83,130.58) ;
	\draw [color={rgb, 255:red, 208; green, 2; blue, 27 }  ,draw opacity=1 ][line width=1.5]    (528.83,130.58) -- (588.83,66.58);
	\draw (547.83,30.58) node    {${\scriptstyle f(x) \ =f_1(x) + f_2(x)}$};
	\draw (531,219) node    {${\scriptstyle \partial _{C} f( 0) \ =\ \{1\}}$};
	\draw (326.83,30.58) node    {${\scriptstyle f_2(x) \ =\ \min\{x,0\}}$};
	\draw (328,217) node    {${\scriptstyle \partial _{C} f_2(0) \ =\ [0,1]}$};
	\draw (122,217) node    {${\scriptstyle\partial _{C} f_1(0) \ =\ [0,1]}$};
	\draw (115.83,30.58) node    {${\scriptstyle f_1(x) \ =\ \max\{x,0\}}$};
	\end{tikzpicture}
\end{tcolorbox}
\vspace{2mm}
	Observe that 
	\[ \del_C f(0) = \{1\} \subsetneq \del_C f_1(0) + \del_C f_2(0) = [0,2]. \]
	The failure of the sum rule is one of the obstacles to computing the Clarke subgradient. Nevertheless, not all is lost, as we still have the following weaker version of the sum rule:
	\[ \del_C (f_1+f_2) \subseteq \del_C f_1 + \del_C f_2; \]
	see~\cite[Proposition 1.12]{Cla75}.
	
	\item[--] {\bf (Tightness).} It is known that if $f$ attains a local minimum at $\bar{\bm x}$, then $\bm{0}\in\del_Cf(\bar{\bm x})$; see~\cite[Proposition 2.3.2]{C90}. By Fact~\ref{fact:clarke}, this is equivalent to $f^\circ(\bar{\bm x},\bm{d}) \ge 0$ for all $\bm{d}\in\R^n$.
	
	However, the Clarke subdifferential may contain stationary points that are not local minima. For instance, consider the function $\R\ni x \mapsto f(x) = -|x|$. It is easy to see that $\del_C f(0) = [-1,1]$. It follows that $\bar{x}=0$ is a stationary point (as $0 \in \del_Cf(0)$). However, the point $\bar{x}=0$ is clearly not a local minimum (in fact, it is a global maximum). Moreover, observe that the corresponding Clarke directional derivatives are $f^\circ(0,1) = f^\circ(0,-1) = 1$, which shows that neither $d=1$ nor $d=-1$ is a descent direction according to Clarke's definition. However, the ordinary directional derivatives exist and are given by $f'(0,1) = f'(0,-1) = -1$. It follows that both $d=1$ and $d=-1$ are descent directions. One may argue that the above example is not persuasive enough, as similar phenomena occur in the smooth case (e.g., $\R \ni x \mapsto f(x) = -x^2$). Hence, let us provide another, perhaps more convincing, example:
	\vspace{3mm}
	\begin{tcolorbox}[breakable]
	 Consider the function
		\[ \R \ni x \mapsto f(x) =\left\{
		\begin{array}{c@{\quad}l}
		x + x^2 \sin (\frac{1}{x}) & \mbox{if } x >0, \\
		x & {\mbox{otherwise}}.
		\end{array} \right.
		\]
\begin{center}		
\begin{tikzpicture}[scale=1.5]
\draw[->](-1.4,0)--(1.4,0)node[left,below,font=\small]{$x$};
\draw[->](0,-1)--(0,1.4)node[right,font=\small]{$f(x)$};
\foreach \x in {-1,0,1}{\draw(\x,0)--(\x,0.05)node[below,outer sep=2pt,font=\small]at(\x,0){\x};}
\foreach \y in {1}{\draw(0,\y)--(0.05,\y)node[left,outer sep=2pt,font=\small]at(0,\y){\y};}
\draw[color={rgb, 255:red, 208; green, 2; blue, 27 }, thick,smooth,domain=0.01:0.7] plot (\x,{\x+\x*\x*sin(1/\x r)});
\draw[color={rgb, 255:red, 208; green, 2; blue, 27 }, thick,smooth,domain=-1:0]plot(\x,\x);
\draw[color={rgb, 255:red, 208; green, 2; blue, 27 },smooth]circle(0.02);
\end{tikzpicture}
\end{center}
	For $x>0$, $f'(x) = 1+2x\sin(\frac{1}{x})-\cos(\frac{1}{x}) $ is bounded on compact sets. Using this and the structure of $f$, it can be shown that $f$ is locally Lipschitz. On one hand, we have $\del _C f(0) = [0,2]$, which means that $\bar{x}=0$ is a stationary point. On the other hand, we have $f'(0,1)=f'(0,-1)=1$. Hence, the point $\bar{x}=0$ is neither a local minimum nor a local maximum.    
	\end{tcolorbox}

Observe that in the above examples, the ordinary directional derivative exists but is strictly smaller than the corresponding Clarke directional derivatives (i.e., $  f'(\bm{x},\bm{d}) < f^\circ(\bm{x},\bm{d})$). This, together with~\eqref{eq:cla-subdiff}, suggests that one may obtain a tighter subdifferential by using other directional derivatives.
\end{itemize}

In view of the aforementioned drawbacks of the Clarke subdifferential, it is natural to ask whether the notion is useful in applications. As it turns out, the Clarke subdifferential can still be a very powerful tool for studying certain sub-classes of locally Lipschitz functions.

\section{Subdifferentially Regular Functions} \label{sec:sd-reg}
In this section, we introduce a representative function class called \emph{subdifferentially regular} functions. The Clarke subdifferential for such functions preserves many of the nice properties of the subdifferential for convex functions. This greatly facilitates the manipulation of such functions in computational procedures. 
\begin{definition}	(\!\!\cite[Definition 2.3.4]{C90}) \label{def:sub-reg}
A locally Lipschitz function $f:\R^n\limto\R$ is \emph{subdifferentially regular} (or simply \emph{regular}) at $\bm{x}\in\R^n$ if for every $\bm{d}\in\R^n$, the ordinary directional derivative \eqref{eq:dir-dev} exists and coincides with the generalized one in \eqref{eq:cla-dir-dev}:
\[f'(\bm{x},\bm{d})  = f^\circ(\bm{x},\bm{d}). \] 
If $f$ is regular at every $\bm{x}\in\R^n$, then we simply say that $f$ is regular.
\end{definition}
As a first example, we note that a convex function $f$ is regular. This follows immediately from~\eqref{eq:cvx-cla}. In this case, we have $\del_Cf = \del f$. Another important example of a regular function is the max function given by $f = \max_{i \in \{1,\ldots,m\}} g_i$, where $g_i:\R^n\limto\R$ ($i=1,\ldots,m$) is smooth; see~\cite[Example 7.28]{RW04}. In particular, this implies that a smooth function is regular. Upon letting $I(\bm{x}) = \{i\in \{1,\ldots,m\}: f(\bm{x}) = g_i(\bm{x})\}$ be the set of indices whose corresponding functions $g_i$ is active at $\bm{x}$, we have $\del_Cf(\bm{x}) =  {\rm conv}\{ \nabla g_i(\bm{x}): i \in I(\bm{x})\}$; see~\cite[Exercise 8.31]{RW04}. We remark that a similar result holds for max functions involving an infinite collection of smooth functions. The interested reader is referred to~\cite[Theorem 10.31]{RW04} for details.

One of the nice properties of regular functions is that they satisfy the following basic calculus rules. 
\begin{fact} \label{prop:cal-rule}
	 (cf.~\cite[Theorem 10.6, Corollary 10.9]{RW04}) 
	 \begin{enumerate}
	 	\item[(a)] {\bf (Chain Rule).} Suppose that $f:\R^n\limto\R$ takes the form $f=g\circ F$, where $g:\R^m \rightarrow \R$ is a locally Lipschitz function and $F:\R^n \rightarrow \R^m$ is a smooth mapping. Given a point $\bm{x} \in \R^n$, if $g$ is regular at $F(\bm{x})$, then $f$ is regular at $\bm{x}$ and  
	 	\[\partial_C f(\bm{x}) = (JF(\bm{x}))^T \partial_C g(F(\bm{x})), \]	
	 	where $JF$ is the Jacobian of $F$; i.e., 
	 	\begin{equation} \label{eq:jac}
	 	JF(\bm{x}) = \left[\frac{\partial f_{i}}{\partial x_{j}}(\bm{x})\right]_{i, j=1}^{m, n} \in \mathbb{R}^{m \times n}.
	 	\end{equation}
        In particular, for a real-valued function $F:\R^n\limto\R$, the gradient of $F$ at $\bm{x}$ is given by $\nabla F(\bm{x}) = (JF(\bm{x}))^T$. 
        
	 \item[(b)] {\bf (Sum Rule).} Suppose that $f  = f_1+f_2+\cdots+f_m$, where $f_i:\R^n \rightarrow \R$ ($i=1,\ldots,m$) are locally Lipschitz functions. Given a point $\bm{x}\in\R^n$, if $f_1,\ldots,f_m$ are regular at $\bm{x}$, then so is $f$ and 
	 \[ \del_C f(\bm{x}) = \del_C f_1(\bm{x}) + \del_C f_2(\bm{x}) + \cdots + \del_C f_m(\bm{x}).\]
	 \end{enumerate}
\end{fact}
We remark that it is possible to develop (possibly weaker) versions of the above calculus rules under weaker assumptions. For instance, a variant of the above chain rule holds in the setting where $g$ is lower semi-continuous\footnote{Recall that a function $f:\R^n\limto\R$ is \emph{lower semi-continuous} if $\liminf_{\bm{y}\limto\bm{x}} f(\bm{y}) = f(\bm{x})$, or equivalently, the epigraph ${\rm epi}(f) = \{(\bm{x},t) \in \R^n\times\R:f({\bm x}) \le t\}$ of $f$ is closed in $\R^n\times\R$; see~\cite[Theorem 1.6]{RW04}.} and $F$ is a locally Lipschitz mapping (and thus not necessarily smooth), while a variant of the above sum rule holds in the setting where $f_1,\ldots,f_m$ are lower semi-continuous. We refer the reader to~\cite[Theorems 10.6 and 10.49]{RW04} for details.

To illustrate the usefulness of the above calculus rules, let us turn our attention to another fundamental class of regular functions, namely weakly convex functions. Such functions have recently received much attention, as they arise in many contemporary signal processing and machine learning applications. We begin with the definition.
\begin{definition}
	A function $f: \R^n \rightarrow\R$ is called \emph{$\rho$-weakly convex} (with $\rho \ge 0$) if the function $\bm{x} \mapsto h(\bm{x}) = f(\bm{x}) + \frac{\rho}{2}\|\bm{x}\|_2^2$ is convex. 
\end{definition}
It is immediate from the definition that a convex function is $0$-weakly convex. As it turns out, weakly convex functions are locally Lipschitz and regular; see~\cite[Propositions 4.4 and 4.5]{V83}. This implies that the basic calculus rules in Fact~\ref{prop:cal-rule} can be applied to weakly convex functions. In particular, we can compute the subdifferential of a weakly convex function $f$ as follows. By definition, the function $\R^n\ni \bm{x}\mapsto h({\bm x})=f(\bm{x}) + \tfrac{\rho}{2}\|\bm{x}\|_2^2$ is convex for some $\rho\ge0$. Using the fact that $\bm{x}\mapsto\tfrac{\rho}{2}\|\bm{x}\|_2^2$ is regular and applying the sum rule, we have $\del_Ch(\bm{x}) = \del_Cf(\bm{x}) + \{\rho\bm{x}\}$. Since $\del_C h$ equals the usual convex subdifferential $\del h$ of $h$, we obtain
\[ \del_C f(\bm{x}) = \del h(\bm{x}) - \{\rho\bm{x}\}. \]

Weakly convex functions are ubiquitous in applications. One prototypical example is the composite function $f=g\circ F$, where $g:\R^m \rightarrow \R$ is convex and Lipschitz continuous on $\R^m$ and $F:\R^n \rightarrow \R^m$ is a smooth map with Lipschitz continuous Jacobian~\cite{DDMC18}. Note that the chain rule in Fact~\ref{prop:cal-rule} yields a formula for $\del_Cf$. Below are some concrete examples of such a composite function that arise in applications.
\begin{itemize}
	\item [--] {\bf (Robust low-rank matrix recovery).} In various signal processing~\cite{davenport2016overview} and machine learning~\cite{srebro2004maximum} applications, a fundamental computational task is to recover a low-rank matrix $\bm{X}^\star \in \R^{n_1\times n_2}$ from a small number of noisy linear measurements of the form
	\[ \bm{y} = \mathcal{A}(\bm{X}^\star) + \bm{s}^\star, \]
	where $\mathcal{A}:\R^{n_1\times n_2}\limto\R^m$ is a known linear operator, $\bm{s}^\star\in\R^m$ is a noise vector, and $\bm{y}\in\R^m$ is the vector of observed values. For simplicity, let us assume that the ground-truth matrix $\bm{X}^\star$ is an $n\times n$ symmetric positive semidefinite matrix of rank $r\ge1$. In the setting where the noise vector represents \emph{outliers} in the measurements, the $\ell_1$-loss function is usually preferred over the $\ell_2$-loss for recovering the ground-truth signal. This gives rise to the following weakly convex formulation for recovering $\bm{X}^\star$~\cite{LZSV20}:
	\[ 
	\min_{\bm{U} \in \R^{n\times r}} f(\bm{U}) = \frac{1}{m} \left\| \bm{y} - \mathcal{A} (\bm{U}\bm{U}^T) \right\|_1.
	\] 
	By applying the chain rule in Fact~\ref{prop:cal-rule}, we can compute
	\[
	\begin{aligned}
	\frac{1}{m} & \big[ (\mathcal{A}^* (\mbox{Sign}(\mathcal{A}(\bm{U}\bm{U}^T)-\bm{y})))^T\bm{U}  \\
	&\,\,\, + \mathcal{A}^* (\mbox{Sign}(\mathcal{A}(\bm{U}\bm{U}^T)-\bm{y}))\bm{U} \big] \subseteq \del_C f(\bm{U}),
	\end{aligned}
	\]
	where $\mathcal{A}^*$ is the adjoint of $\mathcal{A}$; see~\cite{LZSV20}.

	\item [--] {\bf (Robust sign retrieval).} Phase retrieval is a classic inverse problem that arises in areas such as crystallography~\cite{ELB18}, optical imaging~\cite{SEC+15}, and audio signal processing~\cite{Wald17}. Here, let us consider a real-valued version of the problem, in which we are interested in recovering a vector $\bm{x}^\star\in\R^n$ from noisy measurements of the form
	\begin{equation}\label{eq:gen_pr}
	 b_i = (\bm{a}_i^T\bm{x}^\star)^2 + s_i^\star \quad\mbox{for } i=1,\ldots,m, 
	 \end{equation}
	where $\bm{a}_1,\ldots,\bm{a}_m\in\R^n$ are measurement vectors, $\bm{s}^\star\in\R^m$ is the noise vector, and $b_1,\ldots,b_m\in\R$ are the observed values. One approach to tackling this problem is to consider the weakly convex formulation
		\[ \min_{\bm{x} \in \R^n} f(\bm{x})=\frac{1}{m} \sum_{i=1}^{m}\left| (\bm{a_{i}}^T\bm{x})^{2}-b_{i}\right|, \] 
		which aims at handling outliers in the measurements~\cite{DR19,DDP20}. Using the calculus rules in Fact~\ref{prop:cal-rule}, we have
		\[ \frac{2}{m}\sum_{i=1}^m  (\bm{a}_i^T\bm{x}) \cdot {\rm Sign}((\bm{a}_i^T\bm{x})^2-b_i)\cdot \bm{a}_i \subseteq \del_Cf(\bm{x}); \]
		see~\cite{DDP20}.
	
	\item [--] {\bf (Robust blind deconvolution).} 
	The blind deconvolution problem, which is found in diverse fields such as astronomy~\cite{DD1996blind} and image processing~\cite{chan1998total,LWDF11}, aims to recover a pair of signals in two low-dimensional structured spaces from observations of their noisy pairwise convolutions.
	Again, let us focus on a real-valued version of this problem for simplicity. Formally, we consider the task of robustly recovering a pair $(\boldsymbol{w}^\star,\boldsymbol{x}^\star) \in \R^{n_1} \times \R^{n_2}$ from $m$ bilinear measurements:
	\[ b_i = (\boldsymbol{a}_{i}^T\boldsymbol{w}^\star) (\boldsymbol{c}_{i}^T\boldsymbol{x}^\star) + s_i^\star \quad\mbox{for } i=1,\ldots,m, \]
	where $\bm{a}_1,\ldots,\bm{a}_m\in\R^{n_1}$ and $\bm{c}_1,\ldots,\bm{c}_m\in\R^{n_2}$ are measurement vectors,  $b_1,\ldots,b_m\in\R$ are the observed values, and $\bm{s}^\star\in\R^m$ is the noise vector. One non-smooth formulation of the problem reads
\[\min_{\bm{w} \in \R^{n_1},\, \bm{x} \in \R^{n_2}} f(\bm{w},\bm{x})=\frac{1}{m} \sum_{i=1}^{m}\left| (\boldsymbol{a}_{i}^T\boldsymbol{w}) (\boldsymbol{c}_{i}^T\boldsymbol{x}) - b_{i}\right|,\]
in which the $\ell_1$-loss promotes strong recovery and stability guarantees under certain statistical assumptions~\cite{CDDD19}. By invoking the chain rule in Fact~\ref{prop:cal-rule}, we obtain
\[
\begin{aligned}  
\frac{1}{m} &\sum_{i=1}^m\mbox{Sign}  ( (\boldsymbol{a}_{i}^T\boldsymbol{w})(\boldsymbol{c}_{i}^T\boldsymbol{x}) - b_i ) \cdot\\
&\,\,\, \left(  (\boldsymbol{c}_{i}^T\boldsymbol{x})
\left[\begin{array}{l}
\boldsymbol{a}_{i} \\
\boldsymbol{0}
\end{array}\right] + (\boldsymbol{a}_{i}^T\boldsymbol{w}) \left[\begin{array}{l}
\boldsymbol{0} \\
\boldsymbol{c}_{i}
\end{array}\right]\right) \subseteq \del_Cf(\bm{w},\bm{x});
\end{aligned}
\]
see~\cite{CDDD19}.
\end{itemize}

Another illustrative example is given by the family of weakly convex sparse regularizers~\cite{shen2018nonconvex,WCP18}, such as logarithmic sum penalty~\cite{candes2008enhancing}, smoothly clipped absolute deviation (SCAD)~\cite{FL01}, and minimax concave penalty (MCP)~\cite{Zhang10}. These regularizers take the form
\[\R^n \ni \bm{x} \mapsto R(\bm{x}) = \sum_{i=1}^n \phi(|x_i|),\]
where $\phi:\R\limto\R$ is a non-decreasing concave but weakly convex function. Although we cannot apply the chain rule in Fact~\ref{prop:cal-rule} directly, by using the fact that the absolute value function is locally Lipschitz, we can still apply an extended version of the chain rule (see~\cite[Theorem 10.49]{RW04}) to compute an element of the subdifferential of $R$. Let us demonstrate this via the following concrete example:
\vspace{3mm}
\begin{tcolorbox}[breakable]
Let $R:\R^n\limto\R$ be the logarithmic sum penalty function; i.e.,
\[  R(\bm{x})=\sum_{i=1}^{d} \log \left(|x_{i}|+\theta\right), \]
where $\theta>0$ is a smoothing parameter. Consider the following regularized least-squares regression problem:
\[\min _{\bm{x} \in \mathbb{R}^{n}} f(\bm{x}) = \frac{1}{2m} \sum_{i=1}^{m} \left(b_{i}-\bm{a}_{i}^{\top} \bm{x}\right)^2+\lambda R(\bm{x}), \]
where $\lambda\ge0$ is a regularization parameter. Observe that $f$ is regular, as the sum of regular functions is regular; see Fact~\ref{prop:cal-rule}. By the extended chain rule, any vector $\bm{s}\in\R^n$ with
\[ s_i \in \frac{\mbox{Sign}(x_i)}{|x_i|+\theta} \quad\mbox{for } i=1,\ldots,n \]
satisfies $\bm{s} \in \del_C R(\bm{x})$. It is then straightforward to obtain an element of the subdifferential of $f$ via the sum rule in Fact~\ref{prop:cal-rule}.
\end{tcolorbox}

In all the above examples, the ability to explicitly calculate the subdifferential of the weakly convex objective function at hand makes it possible to use simple subgradient methods to minimize the function. Moreover, if the objective function satisfies a regularity condition called \emph{sharpness}, then a suitably initialized subgradient method with properly chosen step sizes will converge at a linear rate to an optimal solution to the problem~\cite{DDMC18} (see also~\cite{LZSV20,DDP20}). We also refer the reader to~\cite{davis2019stochastic}, which discusses stochastic methods for tackling optimization problems involving weakly convex objective functions, and to~\cite{LCD+19}, which develops Riemannian subgradient-type methods for weakly convex optimization over the Stiefel manifold.

Although the class of weakly convex functions provides a powerful modeling tool for applications in signal processing and machine learning, there are still other widely-used functions that do not belong to this class. Here are two examples.
\begin{itemize}
	\item [--] {\bf (Canonical robust sign retrieval).} Besides the squared-amplitude measurement model in~\eqref{eq:gen_pr}, another measurement model of interest for phase retrieval problems is
	\[
	b_i = |\bm{a}_i^T\bm{x}^\star| + s_i^\star \quad\mbox{for } i=1,\ldots,m, 
	\]
	where $\bm{x}^\star \in \R^n$ is the signal to be recovered, $\bm{a}_1,\ldots,\bm{a}_m\in\R^n$ are measurement vectors, $\bm{s}^\star\in\R^m$ is the noise vector, and $b_1,\ldots,b_m\in\R$ are the observed values. Such an amplitude measurement model is used, e.g., in optical wavefront reconstruction; see \cite{LBL02} for details. The corresponding robust phase retrieval problem then takes the form
		\[ \min_{\bm{x} \in \R^n} f(\bm{x})=\frac{1}{m} \sum_{i=1}^{m}\left| |\bm{a_{i}}^T\bm{x}|-b_{i}\right|. \] 
	The function $f$ is not weakly convex as it is not even subdifferentially regular~\cite{ABH19}. 
	
	\item [--] {\bf (Deep Neural Network).} Deep learning is a powerful paradigm in machine learning that allows one to learn a complicated mapping by decomposing it into a series of nested simple mappings, and it has attracted immense interest in various areas of science and engineering~\cite{GBC16}. As an illustration, consider a simple prediction problem, in which one is given $N$ observed feature-label pairs $(\bm{x}_{i},y_i) \in \R^n\times \R$, where $i = 1,\ldots,N$, and the goal is to learn the feature-label relationship. One can model such a relationship using the one-hidden-layer neural network shown in Figure~\ref{fig:one-dnn}. 
\begin{figure}[ht]
	\centering
	\includegraphics[width=0.4\textwidth]{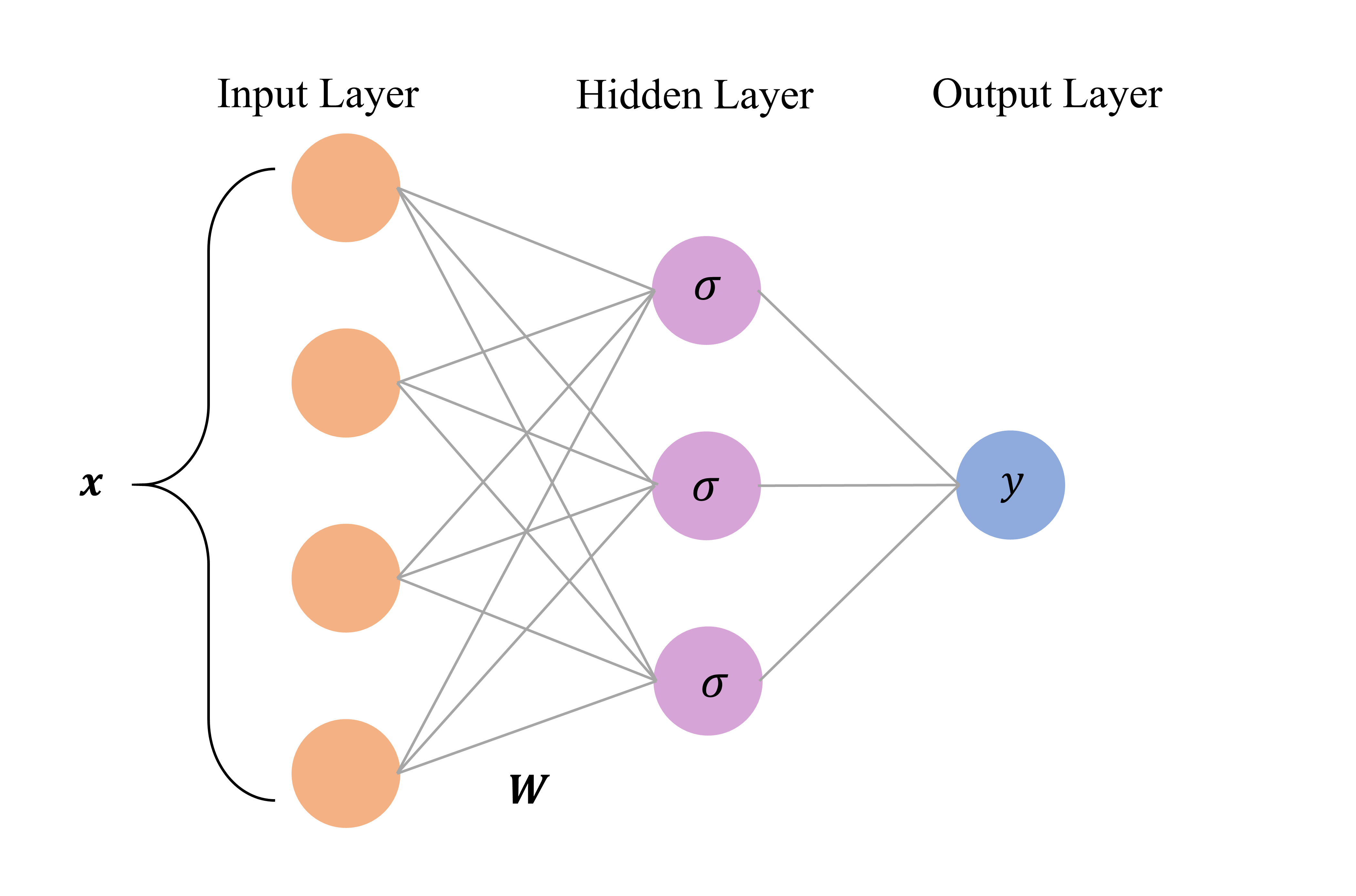}	
	\caption{Illustration of one-hidden-layer neural network.}
	\label{fig:one-dnn}
\end{figure}	
	Here, $\bm{W} = [\bm{w}_1,\ldots,\bm{w}_k]$ is the matrix of weight parameters, where $\bm{w}_j \in \R^n$ denotes the weight with respect to the $j$-th neuron, and $\sigma:\R\limto\R$ is a (typically non-smooth) activation function (e.g., the ReLU function $x\mapsto\max\{x,0\}$). Using the square-loss function, the weights that best model the relationship in the given feature-label pairs can be found by solving the following optimization problem:
	\begin{equation}\label{eq:dnn}
	\min_{\bm{W} \in \R^{n\times k}} f(\bm{W})=\frac{1}{2 N} \sum_{i=1}^{N}\left(\sum_{j=1}^{k} \sigma(\bm{w}_{j}^{T} \bm{x}_{i}) -y_{i}\right)^{2}.
	\end{equation}
	
	Unfortunately, neural networks with non-smooth activation functions typically give rise to objective functions that are not subdifferentially regular~\cite{DDKL20}. For example, consider the instance of Problem~\eqref{eq:dnn} in which $N=n=k=1$, $x_1=y_1=1$, and $\sigma:\R\limto\R$ is the ReLU function (i.e., $\sigma(x)=\max\{x,0\}$). Then, the objective function in~\eqref{eq:dnn} becomes $f(w) = \tfrac{1}{2} (\max\{w,0\}-1)^2$, whose graph is shown below.
	
\begin{center}		
	\begin{tikzpicture}[scale=2.3]
	\draw[->](-1.2,0)--(1.2,0)node[left,below,font=\small]{$w$};
	\draw[->](0,-0.3)--(0,0.8)node[right,font=\small]{$f(w)$};
	\foreach \x in {-1,1}{\draw(\x,0)--(\x,0.05)node[below,outer sep=2pt,font=\small]at(\x,0){\x};}
	\foreach \y in {1/2}{\draw(0,\y)--(0.05,\y)node[left,outer sep=2pt,font=\small]at(0,\y){\y};}
	\draw[color={rgb, 255:red, 208; green, 2; blue, 27 }, thick,smooth,domain=0.001:1] plot (\x,{(1-\x)^2/2});
	\draw[color={rgb, 255:red, 208; green, 2; blue, 27 }, thick,smooth,domain=-1:0]plot(\x,1/2);
	\end{tikzpicture}
\end{center}

It is a simple exercise to show that $f^\circ(0,1)>0>f'(0,1)$. Hence, by Definition \ref{def:sub-reg}, we see that $f$ is not subdifferentially regular. Roughly speaking, the graph of a subdifferentially regular function cannot have ``downward-facing cusps''~\cite{DDKL20}.

\tikzset{every picture/.style={line width=0.75pt}} 
\end{itemize}

In view of the above examples, we are naturally interested in developing other sharper generalized subdifferential concepts that can deal with broader function classes.  

\section{Directionally Differentiable Functions}
As we have seen in Section~\ref{sec:loc-lip}, the Clarke directional derivative $f^\circ$ does not always yield useful information about the descent directions of a function at a given point. For instance, for a directionally differentiable locally Lipschitz function $f$ with directional derivative $f'$, we always have $f' \le f^\circ$ and hence the Clarke subdifferential is in some sense too large; see Fact~\ref{fact:clarke}. We circumvent this problem in Section~\ref{sec:sd-reg} by imposing the assumption $f'=f^\circ$ on the functions we consider, thereby leading us to the class of subdifferentially regular functions. In this section, we present another approach, which begins by constructing subdifferentials that are smaller than the Clarke subdifferential and then trying to refine them so that they possess some of the desirable properties mentioned in Section~\ref{sec:cvx}. One advantage of such an approach is that it allows us to tackle functions that are not necessarily subdifferentially regular.

To begin, consider a directionally differentiable locally Lipschitz function $f:\R^n\limto\R$; i.e., the directional derivative $f'(\bm{x},\bm{d})$ exists for all $\bm{x}\in\R^n$ and $\bm{d} \in \R^n$. Since $f'\le f^\circ$ by Fact~\ref{fact:clarke}, the following set suggests itself as a natural candidate for a subdifferential of $f$:
\begin{equation}\label{eq:fre-sd}
 \widehat{\del}f(\bm{x}) = \left\{ \bm{s} \in \R^n: \bm{s}^T\bm{d} \le f'(\bm{x},\bm{d}) \mbox{ for all } \bm{d} \in \R^n \right\}.
\end{equation}
The set $\widehat{\del}f(\bm{x})$ is known as the \emph{Fr\'{e}chet subdifferential} and its elements the \emph{Fr\'{e}chet subgradients} of $f$ at $\bm{x}$. It is immediate from~\eqref{eq:fre-sd} that $\widehat{\del}f(\bm{x}) \subseteq \del_C f(\bm{x})$ for any $\bm{x}\in\R^n$ (see~\eqref{eq:cla-subdiff}), and that the Fr\'{e}chet subdifferential coincides with the usual convex subdifferential when $f$ is convex (see~\eqref{eq:subdiff-1}).
In fact, the Fr\'{e}chet subdifferential is closely related to the convex subdifferential. Specifically, the former can be obtained by using higher-order minorants in the construction~\eqref{eq:subdiff-2} of the convex subdifferential (see~\cite[Exercises 8.4 and 9.15]{RW04}):
\[ 
\widehat{\del}f(\bm{x}) = \left\{
\bm{s} \in \R^n:
\begin{aligned} & f(\bm{y}) \ge f(\bm{x}) + \bm{s}^T(\bm{y}-\bm{x}) \\
& \quad+ o(\|\bm{y}-\bm{x}\|_2) \mbox{ for all } \bm{y} \in \R^n  
\end{aligned}
\right\}.
\] 
The inequality with the little-oh term in the above expression means that
\[ \liminf_{\bm{y}\limto \bm{x}} \frac{f(\bm{y})-f(\bm{x})-\bm{s}^T(\bm{y}-\bm{x})}{\|\bm{y}-\bm{x}\|_2} \ge 0. \] 
Moreover, observe that for any $\bm{x}\in\R^n$, we have $f'(\bm{x},t\bm{d}) = t\cdot f'(\bm{x},\bm{d})$ for any $\bm{d}\in\R^n$ and $t>0$. Hence, by~\cite[Theorem 8.24]{RW04}, the set $\widehat{\del}f(\bm{x})$ is closed and convex. In addition, since $f'(\bm{x},{\bm d}) < \infty$ for all $\bm{d}\in\R^n$ due to the Lipschitzian property of $f$, the support function of $\widehat{\del}f(\bm{x})$ is given by ${\rm conv}(f'(\bm{x},\cdot))$; i.e.,
\[ {\rm conv}(f'(\bm{x},\cdot))({\bm d}) = \sup_{\bm{s} \in \widehat{\del}f(\bm{x})} \bm{s}^T\bm{d}, \]
where ${\rm conv}(f'(\bm{x},\cdot))$ is the pointwise supremum of all convex functions $g$ satisfying $g(\bm{d}) \le f'(\bm{x},\bm{d})$ for all $\bm{d}\in\R^n$. We refer the interested reader to~\cite{KR03,MNY06} for a detailed treatment of the Fr\'{e}chet subdifferential.

Although the above discussion suggests that the Fr\'{e}chet subdifferential possesses many attractive properties, it is still rather limited. Consider, for instance, the directionally differentiable Lipschitz function $\R\ni x\mapsto f(x)=-|x|$. Then, a simple calculation yields $\widehat{\del}f(0)=\emptyset$. In particular, the Fr\'{e}chet subdifferential can be empty, even at points that could be of interest (in this case, $\bar{x}=0$ is the global maximum). Moreover, by taking a sequence $x^k\searrow0$, we have $-1\in\widehat{\del}f(x^k)$ for all $k$ but $-1\not\in\widehat{\del}f(0)$; i.e., the mapping $\widehat{\del}f$ is not closed. This shows that the Fr\'{e}chet subdifferential is not stable with respect to small perturbations of the point in question, which can cause instabilites in computation. One way of addressing this issue is to ``close'' the mapping $\widehat{\del}f$ by defining the following \emph{limiting subdifferential} of $f$:
\begin{equation} \label{eq:lim-sd}
\del f(\bm{x}) = \left\{
\bm{s} \in \R^n:
\begin{aligned}
& \exists \bm{x}^k \rightarrow \bm{x} ~\mbox{and}~ \bm{s}^k \in \widehat{\del} f(\bm{x}^k)  \\
& \mbox{such that}~ \bm{s}^k \rightarrow \bm{s}
\end{aligned}
\right\}.
\end{equation}
However, such a process can destroy the convexity of the resulting set. Indeed, continuing with the example $f(\cdot)=-|\cdot|$, we have $\del f(0) = \{-1,1\}$. Still, the limiting subdifferential possesses nice properties and is very useful in formulating optimality conditions for non-smooth optimization problems~\cite{KR03,RW04,M13a}. As a first illustration, let us present the following result, which establishes the relationship among the three subdifferentials we have introduced so far, namely the Fr\'{e}chet subdifferential, the limiting subdifferential, and the Clarke subdifferential.
\begin{fact} (cf.~\cite[Theorem 8.6]{RW04}, \cite[Theorem 3.57]{M13a})
	\label{fact:rel}
	For any locally Lipschitz function $f:\R^n\limto\R$ and $\bm{x}\in\R^n$, we have
	\begin{equation} \label{eq:sd-incl}
	\widehat{\partial}f(\bm{x}) \subseteq \partial f(\bm{x}) \subseteq \partial_C f(\bm{x})
	\end{equation}
	and $\del_C f(\bm{x}) = {\rm conv}(\del f(\bm{x}))$. Moreover, if $f$ is subdifferentially regular at $\bm{x}$ (in particular, $f$ is directionally differentiable at $\bm{x}$), then all the above subdifferentials coincide; i.e., $\hat{\partial}f(\bm{x}) = \partial f(\bm{x}) = \partial_C f(\bm{x})$.
\end{fact}

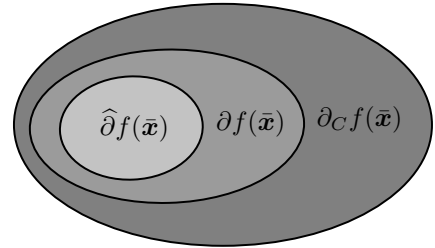
\begin{figure}[h]
	\centering
	\tikzset{every picture/.style={line width=0.75pt}} 
	
	\begin{tikzpicture}[x=0.4pt,y=0.4pt,yscale=-1,xscale=1]
	
	\draw  [fill={rgb, 255:red, 128; green, 128; blue, 128 }  ,fill opacity=1 ] (44,141.7) .. controls (44,78.46) and (132.49,27.2) .. (241.65,27.2) .. controls (350.81,27.2) and (439.3,78.46) .. (439.3,141.7) .. controls (439.3,204.94) and (350.81,256.2) .. (241.65,256.2) .. controls (132.49,256.2) and (44,204.94) .. (44,141.7) -- cycle ;
	\draw  [fill={rgb, 255:red, 155; green, 155; blue, 155 }  ,fill opacity=1 ] (59.05,146.46) .. controls (57.99,106.41) and (115.18,72.42) .. (186.78,70.53) .. controls (258.38,68.64) and (317.28,99.58) .. (318.34,139.63) .. controls (319.39,179.68) and (262.21,213.67) .. (190.61,215.56) .. controls (119.01,217.45) and (60.11,186.51) .. (59.05,146.46) -- cycle ;
	\draw  [fill={rgb, 255:red, 255; green, 255; blue, 255 }  ,fill opacity=0.4 ] (87.53,147.89) .. controls (86.27,120.94) and (115.45,97.69) .. (152.72,95.94) .. controls (189.98,94.2) and (221.21,114.63) .. (222.47,141.57) .. controls (223.73,168.51) and (194.55,191.77) .. (157.28,193.51) .. controls (120.02,195.26) and (88.79,174.83) .. (87.53,147.89) -- cycle ;
	
	\draw (158,142) node    {$ \widehat{\partial} f(\bar{\bm{x}})$};
	\draw (268,140) node    {$\partial  f(\bar{\bm{x}})$};
	\draw (371,137) node    {$\partial _{C} f(\bar{\bm{x}})$};	
	\end{tikzpicture}	
	\caption{Relationship among the various subdifferentials.}
	\label{fig:rel}
\end{figure}

Note that the each of the inclusions in~\eqref{eq:sd-incl} can be strict. Indeed, in our previous example $f(\cdot)=-|\cdot|$, we have $\widehat{\partial} f(0) = \emptyset \subsetneq \partial f(0)  = \{-1,1\} \subsetneq \partial_C f(0) = [-1,1]$. Fact~\ref{fact:rel} reveals that the limiting subdifferential is tighter than the Clarke subdifferential. Moreover, when $f$ is regular, the limiting subdifferential inherits all the properties of the Clarke subdifferential discussed in Sections~\ref{sec:loc-lip} and~\ref{sec:sd-reg}. In particular, since a convex function $f$ is regular, there is no danger of confusion as to the meaning of $\del f$, as the usual convex subdifferential and the limiting subdifferential coincide in this case.

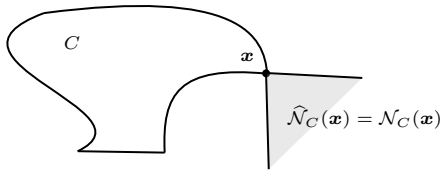
\begin{figure}[h]
	\centering
\tikzset{every picture/.style={line width=0.75pt}}
\begin{tikzpicture}[x=0.75pt,y=0.75pt,yscale=-1,xscale=1]
\draw    (160.34,120) .. controls (158.92,90.13) and (167.41,75.9) .. (212,80.17) ;
\draw    (116.46,119.29) -- (160.34,120) ;
\draw    (99.47,49.58) .. controls (197.85,36.78) and (210.58,60.96) .. (212,80.17) ;
\draw    (99.47,49.58) .. controls (36.49,73.05) and (158.21,99.37) .. (116.46,119.29) ;
\draw  [fill={rgb, 255:red, 0; green, 0; blue, 0 }  ,fill opacity=1 ] (210,80) .. controls (210,79.17) and (210.67,78.5) .. (211.5,78.5) .. controls (212.33,78.5) and (213,79.17) .. (213,80) .. controls (213,80.83) and (212.33,81.5) .. (211.5,81.5) .. controls (210.67,81.5) and (210,80.83) .. (210,80) -- cycle ;
\draw  [fill={rgb, 255:red, 155; green, 155; blue, 155 }  ,fill opacity=0.22 ][line width=0.75]  (212.82,128.38) -- (211.46,79.82) -- (259.99,82.31) ;

\draw (108,60) node [anchor=north west][inner sep=0.75pt]  [font=\scriptsize]  {$C$};
\draw (197,68) node [anchor=north west][inner sep=0.75pt]  [font=\scriptsize]  {$\bm{x}$};
\draw (221,95) node [anchor=north west][inner sep=0.75pt]  [font=\scriptsize]  {$\widehat{\mathcal{N}}_C(\bm{x}) = \mathcal{N}_{C}(\bm{x})$};
\end{tikzpicture}
	\caption{A closed non-convex set with $\widehat{\mathcal{N}}_C(\bm{x}) = \mathcal{N}_C(\bm{x})$.}
\label{fig:n-cone}
\end{figure}

\begin{figure}[h]
	\centering 
\tikzset{every picture/.style={line width=0.75pt}} 

\begin{tikzpicture}[x=0.75pt,y=0.75pt,yscale=-1,xscale=1]

\draw    (139.51,122.99) .. controls (104.49,162.97) and (185.49,167) .. (202.5,135.01) ;
\draw    (202.5,135.01) .. controls (166.48,199) and (271.5,158.03) .. (258.51,128.03) ;
\draw    (139.51,122.99) .. controls (189.53,55) and (247.51,113.02) .. (258.51,128.03) ;
\draw  [fill={rgb, 255:red, 0; green, 0; blue, 0 }  ,fill opacity=1 ] (200.5,136.01) .. controls (200.5,134.9) and (201.4,134.01) .. (202.5,134.01) .. controls (203.61,134.01) and (204.5,134.91) .. (204.5,136.01) .. controls (204.5,137.11) and (203.61,138.01) .. (202.5,138.01) .. controls (201.4,138.01) and (200.5,137.11) .. (200.5,136.01) -- cycle ;
\draw    (92.51,109.97) -- (204.5,136.01) ;
\draw    (203.5,137.01) -- (261.48,193.03) ;
\draw  [dash pattern={on 4.5pt off 4.5pt}]  (133.49,168.98) .. controls (128.95,158.06) and (141.79,140.52) .. (150.17,133.69) ;
\draw [shift={(152.5,131.99)}, rotate = 508.01] [fill={rgb, 255:red, 0; green, 0; blue, 0 }  ][line width=0.08]  [draw opacity=0] (6.25,-3) -- (0,0) -- (6.25,3) -- cycle    ;
\draw  [dash pattern={on 4.5pt off 4.5pt}]  (151.49,174.99) .. controls (181.68,177.75) and (205.19,163.76) .. (215.09,157.53) ;
\draw [shift={(217.5,156.01)}, rotate = 508.01] [fill={rgb, 255:red, 0; green, 0; blue, 0 }  ][line width=0.08]  [draw opacity=0] (6.25,-3) -- (0,0) -- (6.25,3) -- cycle    ;

\draw (120.99,169.98) node [anchor=north west][inner sep=0.75pt]  [font=\tiny,rotate=-0.02]  {$\mathcal{N}_{C}(\bm{x})$};
\draw (198.99,194.01) node [anchor=north west][inner sep=0.75pt]  [font=\tiny,rotate=-0.02]  {$\widehat{\mathcal{N}}_{C}(\bm{x}) =\{\bm{0}\}$};
\draw (167.02,107) node [anchor=north west][inner sep=0.75pt]  [font=\scriptsize,rotate=-0.02]  {$C$};
\draw (203.01,124.01) node [anchor=north west][inner sep=0.75pt]  [font=\scriptsize,rotate=-0.02]  {$\bm{x}$};
\end{tikzpicture}
	\caption{A closed non-convex set with $\widehat{\mathcal{N}}_C(\bm{x}) \subsetneq \mathcal{N}_C(\bm{x})$.}
\label{fig:n-cone-neq}
\end{figure}
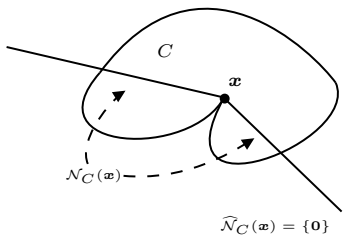

As a further illustration and in preparation for our discussion of optimality conditions of non-smooth optimization problems, let us consider the Fr\'{e}chet and limiting subdifferentials of the indicator function associated with a closed but not necessarily convex set $C\subseteq\R^n$. Recall the definition of the indicator $\mathbb{I}_C$ of $C$ in~\eqref{eq:ind}. Clearly, the indicator needs not be directionally differentiable or locally Lipschitz. Nevertheless, a formal calculation using the definition of the Fr\'{e}chet subdifferential in~\eqref{eq:fre-sd} yields
\begin{equation} \label{eq:fre-sd-ind}
\widehat{\del}\mathbb{I}_C(\bm{x}) = \left\{ 
\bm{s}\in\R^n : 
\begin{aligned}
& \bm{s}^T(\bm{y}-\bm{x}) \le o(\|\bm{y}-\bm{x}\|_2) \\
& \mbox{for all } \bm{y} \in C 
\end{aligned}
\right\} 
\end{equation}
if $\bm{x}\in C$ and $\widehat{\del}\mathbb{I}_C(\bm{x}) = \emptyset$ otherwise. The defining condition of the set on the right-hand side of~\eqref{eq:fre-sd-ind} can also be written as
\[ \limsup_{\bm{y}\limto\bm{x}, \, \bm{y}\in C} \frac{\bm{s}^T(\bm{y}-\bm{x})}{\|\bm{y}-\bm{x}\|_2} \le 0. \]
The formula~\eqref{eq:fre-sd-ind} for $\widehat{\del}\mathbb{I}_C(\bm{x})$ is indeed valid and can be established in a rigorous manner~\cite[Exercise 8.14]{RW04}. The set on the right-hand side of~\eqref{eq:fre-sd-ind} is called the \emph{Fr\'{e}chet normal cone} to $C$ at $\bm{x}$ and is denoted by $\widehat{\mathcal{N}}_C(\bm{x})$; cf.~the discussion following~\eqref{eq:ind-sd}. Now, using~\eqref{eq:lim-sd} and~\eqref{eq:fre-sd-ind}, we can compute the limiting subdifferential of $\mathbb{I}_C$ as
\begin{equation} \label{eq:lim-sd-ind}
\del\mathbb{I}_C(\bm{x}) = \left\{ \bm{s}\in\R^n:
\begin{aligned}
& \exists \bm{x}^k \limto \bm{x} ~\mbox{and}~ \bm{s}^k \in \widehat{\mathcal{N}}_C(\bm{x}^k) \\
& \mbox{such that } \bm{s}^k\limto\bm{s}
\end{aligned}
\right\};
\end{equation}
see~\cite[Definition 6.3 and Exercise 8.14]{RW04}. Following the terminology used above, the set on the right-hand side of~\eqref{eq:lim-sd-ind} is called the \emph{limiting normal cone} of $C$ at $\bm{x}$ and is denoted by $\mathcal{N}_C(\bm{x})$. Figures~\ref{fig:n-cone} and~\ref{fig:n-cone-neq} show the  Fr\'{e}chet and limiting normal cones of two closed non-convex sets. It is worth noting that the two normal cones do not always coincide; see Figure~\ref{fig:n-cone-neq}, where $\widehat{\mathcal{N}}_C(\bm{x})$ consists of the zero vector only and $\mathcal{N}_C(\bm{x})$ consists of the two rays emanating from $\bm{x}$. In general, we always have $\widehat{\mathcal{N}}_C(\bm{x}) \subseteq \mathcal{N}_C(\bm{x})$~\cite[Proposition 6.5]{RW04}.

\subsection{Concepts of Stationarity}
Armed with the above development, we are now ready to address our primary goal of this paper, which is to introduce and compare different stationarity concepts for non-convex non-smooth optimization problems. To begin, consider Problem~\eqref{eq:con-opt}, where $g:\R^n\limto\R$ is a directionally differentiable locally Lipschitz function and $C\subseteq\R^n$ is a closed set. We say that $\bar{\bm x}\in\R^n$ is a \emph{directional stationary} (resp.~\emph{limiting stationary} and~\emph{Clarke stationary}) point of Problem~\eqref{eq:con-opt} if $\bm{0} \in \widehat{\del}(g+\mathbb{I}_C)(\bm{x})$ (resp.~$\bm{0}\in\del(g+\mathbb{I}_C)(\bm{x})$ and~$\bm{0}\in\del_C(g+\mathbb{I}_C)(\bm{x})$). The following result gives a necessary condition for local optimality of a feasible solution to Problem~\eqref{eq:con-opt}:
\begin{fact} (cf.~\cite[Theorems 8.15 and 10.1, Corollary 6.29]{RW04})
\label{fact:fre-opt-cond}
If $\bar{\bm x}$ is a local minimum of~\eqref{eq:con-opt}, then $\bar{\bm x}$ is a directional stationary (d-stationary) point of~\eqref{eq:con-opt}. If in addition $g$ and $\mathbb{I}_C$ are regular at $\bar{\bm x}$, then 
\[ f'(\bar{\bm x},\bm{d}) \ge 0 \quad\mbox{for all } \bm{d} \in \mathcal{N}_C^\circ(\bm{x}), \]
where 
\[ \mathcal{N}_C^\circ(\bm{x}) = \left\{ \bm{d}\in\R^n : \bm{s}^T\bm{d} \le 0 \mbox{ for all } \bm{s}\in\mathcal{N}_C(\bm{x}) \right\} \]
is called the \emph{polar} of $\mathcal{N}_C(\bm{x})$.
\end{fact}

Note that if $\bar{\bm x}$ is a d-stationary point of~\eqref{eq:con-opt}, then by Facts~\ref{fact:rel} and~\ref{fact:fre-opt-cond} it is also a limiting stationary (l-stationary) and Clarke stationary (C-stationary) point of~\eqref{eq:con-opt}. In particular, we have the following implications:
\[ \mbox{d-stationarity} \Longrightarrow \mbox{l-stationarity}\Longrightarrow \mbox{C-stationarity}. \] 
We now give two examples to show that the reverse implications need not hold in general; see~\cite{CPS18}.
\vspace{2mm}
\begin{tcolorbox}[breakable]

	\tikzset{every picture/.style={line width=0.75pt}} 
	
	\begin{tikzpicture}[x=0.5pt,y=0.5pt,yscale=-1,xscale=1]
	
	\draw [color={rgb, 255:red, 73; green, 135; blue, 206 }  ,draw opacity=1 ][line width=1.5]  (53,162.1) -- (216.5,162.1)(134.22,93.89) -- (134.22,225.58) (209.5,157.1) -- (216.5,162.1) -- (209.5,167.1) (129.22,100.89) -- (134.22,93.89) -- (139.22,100.89) (165.22,157.1) -- (165.22,167.1)(196.22,157.1) -- (196.22,167.1)(103.22,157.1) -- (103.22,167.1)(72.22,157.1) -- (72.22,167.1)(129.22,131.1) -- (139.22,131.1)(129.22,193.1) -- (139.22,193.1) ;
	\draw   ;
	\draw [color={rgb, 255:red, 208; green, 2; blue, 27 }  ,draw opacity=1 ][line width=1.5]    (79.5,212.89) -- (134.22,162.1) ;
	\draw [shift={(134.22,162.1)}, rotate = 317.13] [color={rgb, 255:red, 208; green, 2; blue, 27 }  ,draw opacity=1 ][fill={rgb, 255:red, 208; green, 2; blue, 27 }  ,fill opacity=1 ][line width=1.5]      (0, 0) circle [x radius= 4.36, y radius= 4.36]   ;
	
	\draw [color={rgb, 255:red, 208; green, 2; blue, 27 }  ,draw opacity=1 ][line width=1.5]    (134.22,162.1) -- (164.88,192.4) ;

	\draw [color={rgb, 255:red, 208; green, 2; blue, 27 }  ,draw opacity=1 ][line width=1.5]    (164.88,192.4) -- (213.5,147.89) ;
	
	\draw [shift={(164.88,192.4)}, rotate = 317.52] [color={rgb, 255:red, 208; green, 2; blue, 27 }  ,draw opacity=1 ][fill={rgb, 255:red, 208; green, 2; blue, 27 }  ,fill opacity=1 ][line width=1.5]      (0, 0) circle [x radius= 4.36, y radius= 4.36]   ;
	
	\draw [color={rgb, 255:red, 73; green, 135; blue, 206 }  ,draw opacity=1 ][line width=1.5]  (291.5,163.1) -- (484.5,163.1)(387.37,94.89) -- (387.37,226.58) (477.5,158.1) -- (484.5,163.1) -- (477.5,168.1) (382.37,101.89) -- (387.37,94.89) -- (392.37,101.89) (418.37,158.1) -- (418.37,168.1)(449.37,158.1) -- (449.37,168.1)(356.37,158.1) -- (356.37,168.1)(325.37,158.1) -- (325.37,168.1)(382.37,132.1) -- (392.37,132.1)(382.37,194.1) -- (392.37,194.1) ;
	\draw   ;
	\draw [color={rgb, 255:red, 208; green, 2; blue, 27 }  ,draw opacity=1 ][line width=1.5]    (324.65,163.89) -- (387.37,163.1) ;
	\draw [shift={(387.37,163.1)}, rotate = 359.28] [color={rgb, 255:red, 208; green, 2; blue, 27 }  ,draw opacity=1 ][fill={rgb, 255:red, 208; green, 2; blue, 27 }  ,fill opacity=1 ][line width=1.5]      (0, 0) circle [x radius= 4.36, y radius= 4.36]   ;
	
	\draw [color={rgb, 255:red, 208; green, 2; blue, 27 }  ,draw opacity=1 ][line width=1.5]    (387.37,163.1) -- (449.65,215.89) ;

	\draw [color={rgb, 255:red, 208; green, 2; blue, 27 }  ,draw opacity=1 ][line width=1.5]    (324.65,163.89) -- (276.65,116.89) ;
	
	\draw [shift={(324.65,163.89)}, rotate = 224.4] [color={rgb, 255:red, 208; green, 2; blue, 27 }  ,draw opacity=1 ][fill={rgb, 255:red, 208; green, 2; blue, 27 }  ,fill opacity=1 ][line width=1.5]      (0, 0) circle [x radius= 4.36, y radius= 4.36]   ;

	\draw (323,178) node  [font=\footnotesize,color={rgb, 255:red, 74; green, 144; blue, 226 }  ,opacity=1 ]  {$-1$};
	\draw (378,152) node  [font=\footnotesize,color={rgb, 255:red, 74; green, 144; blue, 226 }  ,opacity=1 ]  {$0$};
	\draw (382.83,71.58) node    {${\scriptstyle f_{2}(x) \ =\ \max\{ -x-1, \min\{ -x,0 \} \}}$};
	\draw (127,152) node  [font=\footnotesize,color={rgb, 255:red, 74; green, 144; blue, 226 }  ,opacity=1 ]  {$0$};
	\draw (165,173) node  [font=\footnotesize,color={rgb, 255:red, 74; green, 144; blue, 226 }  ,opacity=1 ]  {$0.5$};
	\draw (144.83,70.58) node    {${\scriptstyle f_{1}(x) \ =\ \max\{ -|x|,x-1 \}}$};
	\end{tikzpicture}
		\begin{itemize}
		\item[--] For the univariate function $f_1:\R\limto\R$, we have $\del_C f_1(0) = [-1,1]$ and $\del f_1(0) = \{-1,1\}$. It follows that the point $\bar{x}=0$ is C-stationary but fails to be l-stationary. The unique l-stationary point is $x^\star=0.5$ and is also a local minimum. 
		
		\item [--] For the univariate function $f_2:\R\limto\R$, we have $\del f_2(0) = \{-1,0\}$ and $\widehat{\del} f_2(0) = \emptyset$. It follows that the point $\bar{x}=0$ is l-stationary but not d-stationary. The unique d-stationary point is $x^\star=-1$ and is also a local minimum.
	\end{itemize}	
\end{tcolorbox}
\vspace{2mm}

The above discussion suggests that among the three notions of stationarity, d-stationarity is the sharpest. However, the development of algorithms for computing a d-stationary point of the non-convex non-smooth optimization problem~\eqref{eq:con-opt} is still in the infancy stage. We will briefly discuss a recent effort in this direction in the next sub-section and refer the reader to~\cite{PRA17,CPS18} for further reading. By contrast, under the assumption that $g+\mathbb{I}_C$ satisfies the so-called \emph{Kurdyka-{\L}ojasiewicz} property, various algorithms will produce iterates that are provably convergent to a limiting stationary point of~\eqref{eq:con-opt}; see, e.g.,~\cite{ABS13}.

\subsection{Application: Least Squares Piecewise Affine Regression}
In this sub-section, we discuss a representative application called \emph{Least Squares Piecewise Affine Regression}, in which the objective function is piecewise linear-quadratic (PLQ) and hence directionally differentiable (see~\cite[Proposition 10.21]{RW04}).  Specifically, the objective function takes the form
\begin{equation}\label{eq:LSPAR}
\min\limits_{\bm{W} \in C} 
f(\bm{W}) = \frac{1}{2N} \sum_{s=1}^N \left( y_s - \max\limits_{1\leq i \leq k}\bm{w}_i^T\bm{x}_s \right)^2,
\end{equation}
where $\bm{W}=[\bm{w}_1,\ldots,\bm{w}_k] \in \R^{n\times k}$ is the matrix of decision variables and $C\subseteq\R^{n\times k}$ is the feasible set. By setting $h_s(u) = (y_s-u)^2$ (the square loss) and $g_s(\bm{W}) = \max_{1\le i\le k} \bm{w}_i^T\bm{x}_s$ (a piecewise affine function), we can write the above problem in the following compact form:
\[ 
\min\limits_{\bm{W}\in C} f(\bm{W}) =  \frac{1}{2N}\sum\limits_{s=1}^N h_s(g_s(\bm{\theta})).
\] 
The above problem can be used to model the one-layer neural network with the ReLU activation function, in which $k=1$, $C=\R^n$, and $g_s$ takes the simple form $g_s(\bm{w}) = \max\{\bm{w}^T\bm{x}_s,0\}$; cf.~\eqref{eq:dnn}. Our interest in Problem~\eqref{eq:LSPAR} stems from the following:
\begin{fact} (cf.~\cite[Proposition 16]{cui2018finite}) 
The least squares piecewise affine regression problem \eqref{eq:LSPAR} possesses the following properties:
\begin{enumerate}
	\item[(a)] It attains a finite global minimum value.
	\item[(b)] The set of d-stationary points is finite. 
	\item[(c)] Every d-stationary point is a local minimizer.
\end{enumerate}
\end{fact}
The above result provides further evidence that the notion of d-stationarity is in some sense the sharpest, as every d-stationary point of Problem~\eqref{eq:LSPAR} is a local minimum. In view of this, it is natural to ask whether we can propose an iterative algorithm to find such points. In \cite{CPS18} the authors proposed a non-monotone majorized-minimization (MM) algorithm with a semi-smooth Newton method as its inner solver to find a d-stationary point of a class of so-called \emph{composite difference-convex-piecewise} optimization problems, of which Problem~\eqref{eq:LSPAR} is an instance. They also showed that the MM algorithm will converge to a d-stationary point of such problems under mild conditions (which are satisfied by~\eqref{eq:LSPAR}). One of the motivations for introducing such an algorithm is that it is not known whether the basic chain rule holds for the objective function in~\eqref{eq:LSPAR}. For the purpose of experimentation, let us pretend the basic chain rule holds and use it to compute a pseudo-subgradient (actually back-propagation in deep learning) of the objective function:
\begin{equation*}
\frac{\widetilde{\del} f(\bm{w}_i)}{\del \bm{w}_i} = \frac{1}{2} \sum_{s=1}^N \left( \max\limits_{1\leq i \leq k}\bm{w}_i^T\bm{x}_s-y_s \right)x_s \mathbb{I}_{\left\{ i \in \arg\max\limits_{i}\bm{w}_i^T\bm{x}_s \right\}}.
\end{equation*}
Then, we can try using the subgradient method with such a pseudo-subgradient to tackle Problem~\eqref{eq:LSPAR}. However, such an approach does not quite work empirically. 
\begin{figure}[ht]
	\centering
	\includegraphics[scale=0.5]{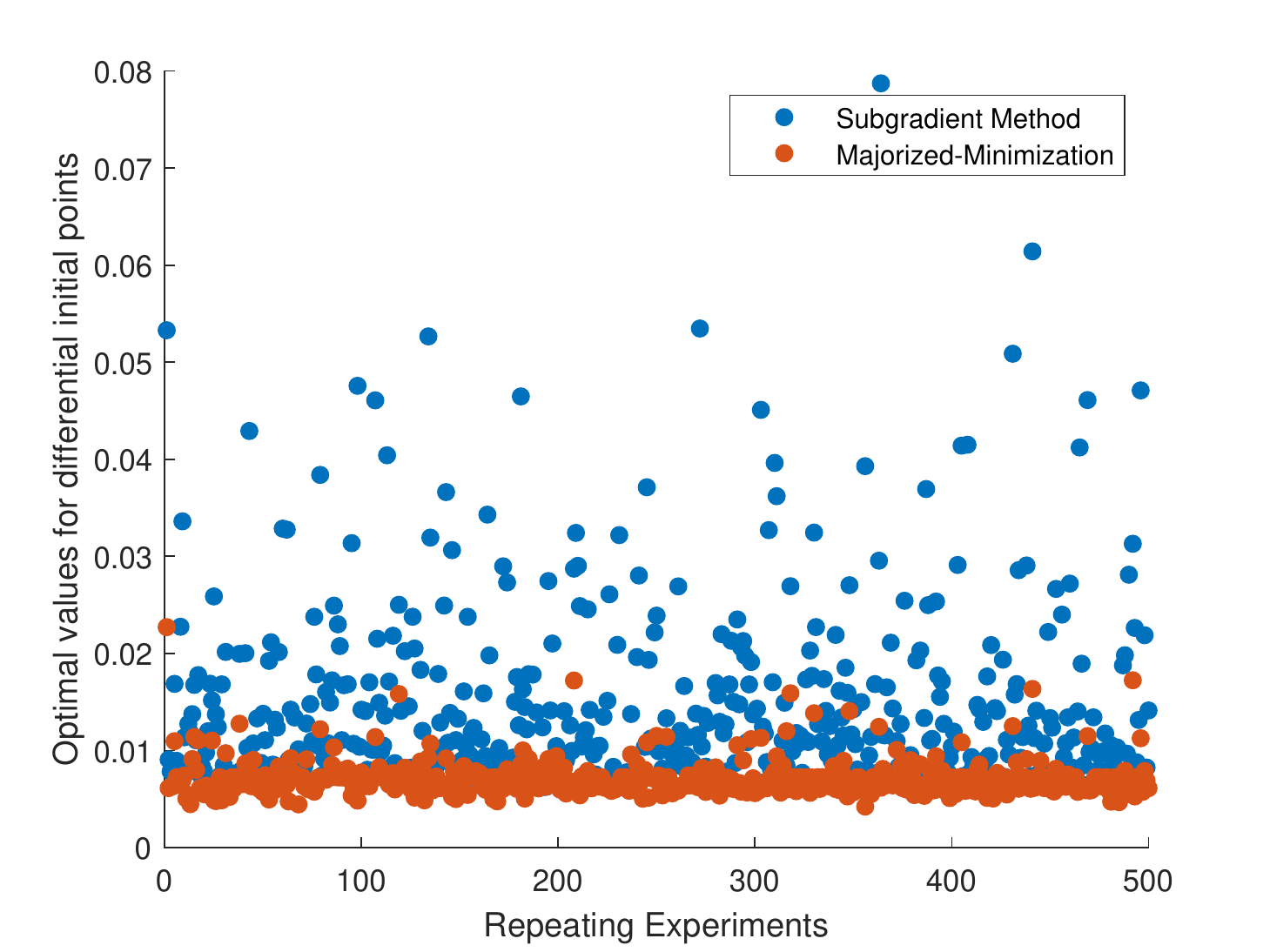}
	\caption{Objective values computed by the MM and subgradient algorithms, $N = 10$.}
	\label{fig:label}
\end{figure}
Indeed, let us follow the experimental setup in~\cite{CPS18} and consider the 2-dimensional convex piecewise linear model
\[y=\max \left\{x_{1}+x_{2}, x_{1}-x_{2},-2 x_{1}+x_{2},-2 x_{1}-x_{2}\right\}+\varepsilon \]
with different sample sizes $N = 10,50,100$. We test the MM and subgradient algorithms on synthetic data. Using the same initial points for the two algorithms, all the experiment results reported here were collected over 500 independent trials over random seeds. 
\begin{figure}[ht]
	\centering
	\includegraphics[scale=0.5]{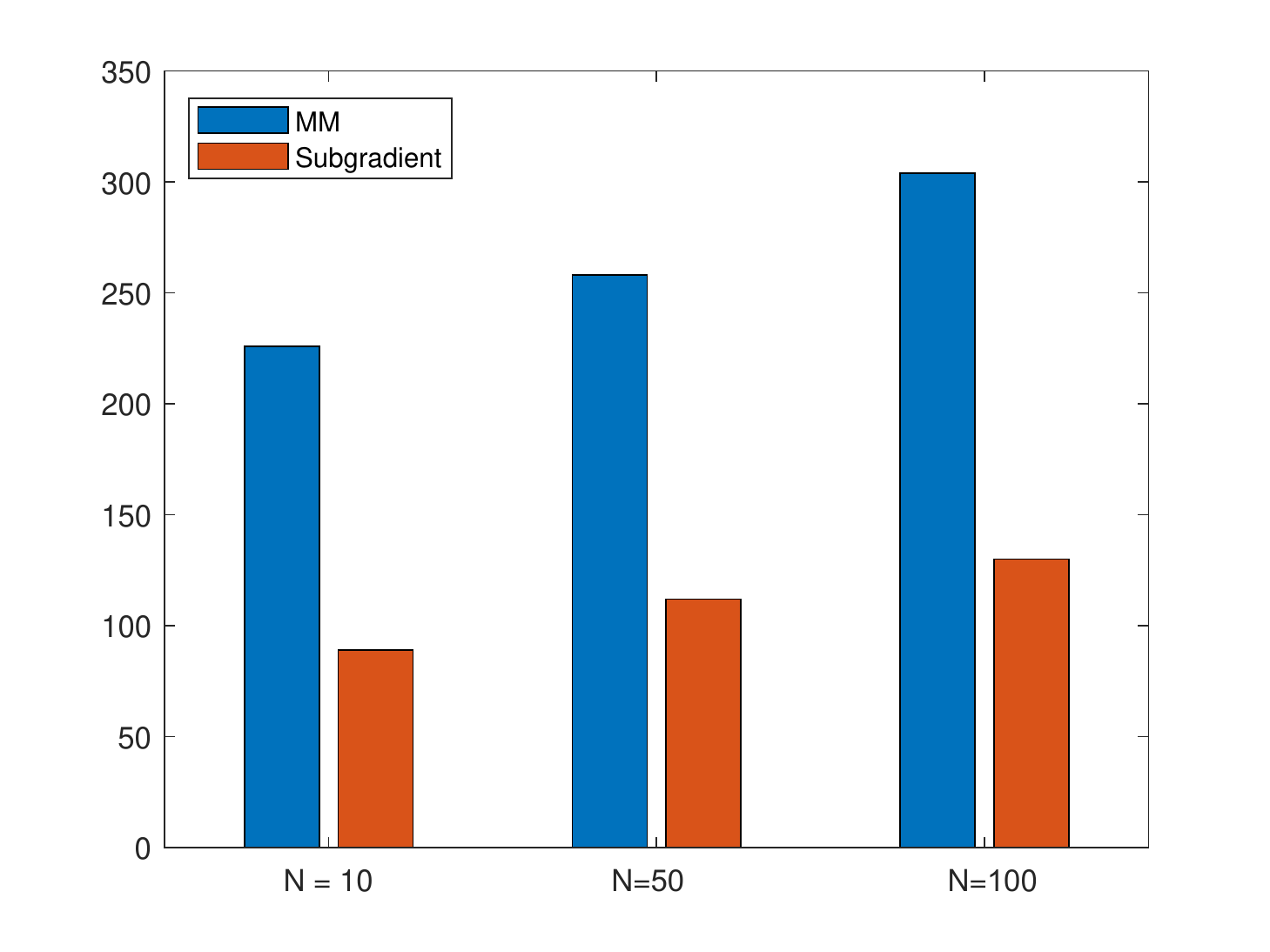}
	\caption{Number of initial points that lead to the smallest objective values.}
	\label{fig:2}
\end{figure}
From Figures~\ref{fig:label} and~\ref{fig:2}, we observe that there is an apparent gap between these two algorithms. In particular, the figures show that the subgradient algorithm reaches many limit points that are unsatisfactory.  Nevertheless, the MM algorithm can be rather slow. As a future work, it would be interesting to design practically efficient first-order algorithms that can provably return a d-stationary point for this application, and more generally, for other signal processing, machine learning, and statistical applications; see, e.g.,~\cite{APX17,PRA17,cui2018finite,CPS18,NPR19} and the references therein.

\section{Conclusion}
In this article, we elucidated the constructions of various subdifferentials for several important sub-classes of non-smooth functions and discussed their corresponding stationarity concepts. We also showcased several representative examples and applications to illustrate the differences among various constructions. We hope that this introductory article will serve as a good starting point for readers who would like to utilize the mathematical tools from non-smooth analysis in the design and analysis of iterative methods for non-smooth optimization problems.

%


%

%
%
%
%
%
%

\ifCLASSOPTIONcaptionsoff
  \newpage
\fi



%

%


\bibliography{sdpbib}
\bibliographystyle{IEEEtran}

\end{document}